\documentclass[a4paper, 10pt, DIV10, headinclude=false, footinclude=false, leqno]{scrartcl}

\usepackage{amsmath, amssymb, amsfonts, amsthm}
\usepackage[utf8]{inputenc}
\usepackage[english]{babel}
\usepackage{url}
\usepackage{mathtools}
\usepackage{esint}
\usepackage{pgfplots}
\pgfplotsset{compat=newest}
\usepgfplotslibrary{groupplots}
\usepgfplotslibrary{dateplot}
\usepackage{subcaption} 
\usepackage[square, numbers]{natbib}

\newcommand{\diff}[2]{\frac{{\partial #1}}{{\partial #2}} }

\usepackage{tikz}
\usetikzlibrary{spy}
\usepackage{verbatim}
\usetikzlibrary{arrows,positioning,calc}
\usepackage{ifthen}
\usepackage{xintexpr}
\tikzset{
    >=stealth',
    punkt/.style={
           rectangle,
           draw=black, thick,
           text width=6.5em,
           minimum height=12em,
           minimum width=12em,
           fill = {rgb:red,108;green,152;blue,221},
		   },
    pil/.style={
           ->,
           shorten <=5pt,
           shorten >=5pt, anchor=north}
}
\usepackage{pgfplots}

\usetikzlibrary{plotmarks}

\usepackage{hyperref}

\numberwithin{figure}{section}
\numberwithin{table}{section}
\numberwithin{equation}{section}

\newenvironment{abstr}[1]{ \vspace{.05in}\footnotesize
       \parindent .2in
         {\upshape\bfseries #1. }\ignorespaces}{\par\vspace{.1in}}
\newenvironment{Abstract}{\begin{abstr}{Abstract}}{\end{abstr}}
\newenvironment{keywords}{\begin{abstr}{Key words}}{\end{abstr}}

\usepackage{thmtools}

\declaretheoremstyle[%
	headindent = \parindent,
	headfont=\sffamily\bfseries,
	notefont=\normalfont\sffamily,
	bodyfont=\normalfont,
	headformat=\NUMBER\ \NAME\NOTE,
	headpunct={},
	postheadspace=1ex,
	spaceabove=0pt,spacebelow=0pt,]%
{mainstyle}
\declaretheoremstyle[%
	headfont=\bfseries\scshape,
	bodyfont=\normalfont,
	headpunct=:,
	postheadspace=1ex,
	spacebelow=12pt,spaceabove=2pt,
	qed=\qedsymbol]%
{beweise}

\declaretheorem[name=Definition,parent=section,style=mainstyle]{definition}
\declaretheorem[name=Remark,sharenumber=definition,style=mainstyle]{bemerkung}

\declaretheorem[name=Lemma,sharenumber=definition,style=mainstyle]{lemma}

\declaretheorem[name=Theorem,sharenumber=definition,style=mainstyle]{theorem}

\declaretheorem[name=Proof,numbered=no,style=beweise]{beweis}

\DeclarePairedDelimiterX\mengenA[1]{\lbrace}{\rbrace}{#1}
\DeclarePairedDelimiterX\mengenB[2]{\lbrace}{\rbrace}{#1\, \delimsize\vert \, #2}

\makeatletter
\newcommand{\set}[2][\relax]{
\ifx#1\relax \ensuremath{
\mengenA*{#2}}
\else \ensuremath{%
  \mengenB*{#1}{#2}}
\fi}
\makeatother

\DeclareRobustCommand{\minwidthbox}[2]{%
  \ifmmode
    \expandafter\mathmakebox
  \else
    \expandafter\makebox
  \fi
  [\ifdim#2<\width\width\else#2\fi]{#1}%
}
\DeclarePairedDelimiter{\abs}{|}{|}
\DeclarePairedDelimiterX\skal[2]{(}{)}{#1\,,\,#2}
\DeclarePairedDelimiter{\norm}{\lVert}{\rVert}
\newcommand{\vertiii}[1]{{\vert\kern-0.3ex\vert\kern-0.3ex\vert #1
    \vert\kern-0.3ex\vert\kern-0.3ex\vert}}  %energynorm

\DeclareMathOperator{\id}{Id} %identische Abbildung
\DeclareMathOperator{\esssup}{ess \, sup}

\DeclareMathOperator{\essinf}{ess \, inf}

\newcommand{\dd}{{\text d}}

\newcommand{\QQ}{{\mathcal{Q}}}
\newcommand{\RR}{{\mathcal{R}}}
\newcommand{\TT}{{\mathcal{T}}}
\newcommand{\IH}{\mathcal{I}_H}
\DeclareRobustCommand{\rchi}{{\mathpalette\irchi\relax}}
\newcommand{\irchi}[2]{\raisebox{\depth}{$#1\chi$}} % inner command, used by \rchi

\usepackage{amsmath}
\usepackage[ruled,linesnumbered]{algorithm2e}
\let\oldnl\nl
\newcommand{\nonl}{\renewcommand{\nl}{\let\nl\oldnl}}

\allowdisplaybreaks[4]

\begin{document}

\title{Numerical upscaling of perturbed diffusion problems}

\author{Fredrik Hellman\footnotemark[1] \and Tim Keil\footnotemark[2] \and Axel Målqvist\footnotemark[1]}
\date{\today}
% \todo{Remove table of contents in the final version}
% \tableofcontents
\newpage
\maketitle

\renewcommand{\thefootnote}{\fnsymbol{footnote}}
%\footnotetext[3]{?}
\footnotetext[2]{Institute for Computational and Applied Mathematics, Westf\"alische Wilhelms-Uni\-ver\-si\-t\"at M\"unster, Einsteinstr. 62, D-48149 M\"unster, Germany. Funded by Deutsche Forschungsgemeinschaft (DFG, German Research Foundation) under Germany
	’s Excellence Strategy – EXC 2044 – 390685587, Mathematics Münster and by the DFG under contract SCHI 1493/1-1.}
\footnotetext[1]{Department of Mathematical Sciences, Chalmers University of Technology and University of Gothenburg SE-421 96 Göteborg, Sweden. The first and third authors work was supported by the Swedish Research Council under Grant 2015-04964 and the G\"{o}ran Gustafsson Foundation for research in natural sciences and medicine.}
\renewcommand{\thefootnote}{\arabic{footnote}}

\begin{Abstract}
In this paper we study elliptic partial differential equations with rapidly varying diffusion coefficient that can be represented as a perturbation of a reference coefficient. We develop a numerical method for efficiently solving multiple perturbed problems by reusing local computations performed with the reference coefficient. The proposed method is based on the Petrov--Galerkin Localized Orthogonal Decomposition (PG-LOD) which allows for straightforward parallelization with low communcation overhead and memory consumption. We focus on two types of perturbations: local defects which we treat by recomputation of multiscale shape functions and global mappings of a reference coefficient for which we apply the domain mapping method. We analyze the proposed method for these problem classes and present several numerical examples.
\end{Abstract}


\section{Introduction}
\label{sec:introduction}
Manufactured heterogeneous materials, such as composites with tailored properties, are crucial tools in engineering. The challenge of performing accurate computer simulations involving such materials have driven the development of multiscale methods over decades \cite{hug,EEnq,MsFEM,MP14,Gamblets}. Multiscale methods have turned out to be successful in computing coarse-scale representations of the solutions to such problems. However, when the heterogeneous data is perturbed it is not obvious how multiscale methods can be adapted. Understanding the effect of perturbations is important since manufactured materials, in general, will not be perfect. Manufacturing tolerances and faults lead to perturbations in the material distribution. There are also other problems, such as time stepping with time dependent diffusion coefficient, optimization of material distribution, and non-linear diffusion problems, which call for iterative procedures where the data in the current iterate can be seen as a perturbation of the data in the previous.

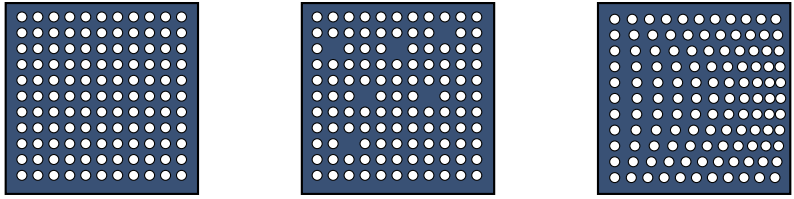
\begin{figure}
\begin{center}
	\begin{tikzpicture}[scale=0.6, every node/.style={transform shape}, node distance=1cm]
	 %nodes
	 \node[] (market) {};
	 \node[above=of market] (dummy) {};
	 \node[punkt, right=of dummy] (t) {};
	 \node[punkt, left=of dummy] (g) {};
	 \node[left=of g] (sec_dummy) {};
	 \node[punkt, left=of sec_dummy] (a) {};
	 \node[below=of sec_dummy](sec_market) {};
	 \foreach \x in {-11.45,-11.1,-10.75,...,-7.6}{
	 \foreach \y in {-0.45,-0.1,0.25,...,3.4}{
	 		\draw[fill=white] (-0.045+\x,\y) circle [radius=0.11, color=white];
	 	}
	}
	\foreach \x in {1.55,1.9,2.25,...,5.4}{
	\foreach \y in {0.1,0.45,...,3.8}{
			\draw[fill=white] (-0.045+\x+4* 1.7 * \y / 3.85 *\x/3.85-4* 1.7 * \y / 3.85 * 1.55/3.85 -4* 1.7 * \y / 3.85 * \x^2/3.85^2 +4* 1.7 * \y / 3.85 * 2* 1.55 * \x / 3.85^2 -4* 1.7 * \y / 3.85 * 1.55^2/3.85^2 - 4* 1.7 * \y^2 / 3.85^2 *\x/3.85 + 4* 1.7 * \y^2 / 3.85^2 *1.55/3.85 + 4* 1.7 * \y^2 / 3.85^2 *\x^2/3.85^2 -4* 1.7 * \y^2 / 3.85^2 * 2* 1.55 * \x / 3.85^2 + 4* 1.7 * \y^2 / 3.85^2 *1.55^2/3.85^2 ,\y - 0.6)  circle [radius=0.11, color=white];
		}
	}

	%TODO : REPLACE THIS BY USING IF STATEMENTS

		\foreach \x in {-5.}{
		\foreach \y in {-0.45,-0.1,0.25,...,3.4}{
			\draw[fill=white] (-0.02+\x,\y) circle [radius=0.11, color=white];
		}
	}
	
			\foreach \x in {-4.65}{
	\foreach \y in {-0.45,-0.1,...,2.}{
		\draw[fill=white] (-0.02+\x,\y) circle [radius=0.11, color=white];
	}
}

			\foreach \x in {-4.65}{
	\foreach \y in {2.7,3.05}{
		\draw[fill=white] (-0.02+\x,\y) circle [radius=0.11, color=white];
	}
}

	\foreach \x in {-4.3}{
	\foreach \y in {-0.45,-0.1}{
			\draw[fill=white] (-0.02+\x,\y) circle [radius=0.11, color=white];
		}
	}

	\foreach \x in {-4.3}{
	\foreach \y in {0.6,0.95,...,3.05}{
			\draw[fill=white] (-0.02+\x,\y) circle [radius=0.11, color=white];
		}
	}

	\foreach \x in {-3.95,-3.6,...,-1.4}{
	\foreach \y in {-0.45,-0.1,0.25,0.6,0.95}{
			\draw[fill=white] (-0.02+\x,\y) circle [radius=0.11, color=white];
		}
	}
	
	\foreach \x in {-3.95}{
	\foreach \y in {1.65,2.,...,3.4}{
			\draw[fill=white] (-0.02+\x,\y) circle [radius=0.11, color=white];
		}
	}
		
	\foreach \x in {-3.60}{
	\foreach \y in {1.3,1.65,2.,...,3.4}{
			\draw[fill=white] (-0.02+\x,\y) circle [radius=0.11, color=white];
		}
	}
		
	\foreach \x in {-3.25}{
	\foreach \y in {1.3,1.65,2.,2.7,3.05}{
			\draw[fill=white] (-0.02+\x,\y) circle [radius=0.11, color=white];
		}
	}
	
	\foreach \x in {-2.9}{
	\foreach \y in {1.3,1.65,2.,...,3.4}{
			\draw[fill=white] (-0.02+\x,\y) circle [radius=0.11, color=white];
		}
	}

	\foreach \x in {-2.55}{
	\foreach \y in {1.65,2.,2.35,2.7,3.05}{
			\draw[fill=white] (-0.02+\x,\y) circle [radius=0.11, color=white];
		}
	}
	
	\foreach \x in {-2.2}{
	\foreach \y in {1.3,1.65,2.,2.35,3.05}{
			\draw[fill=white] (-0.02+\x,\y) circle [radius=0.11, color=white];
		}
	}
			
	\foreach \x in {-1.85}{
	\foreach \y in {1.3,1.65,2.,...,3.4}{
			\draw[fill=white] (-0.02+\x,\y) circle [radius=0.11, color=white];
		}
	}
		
	\foreach \x in {-1.5}{
	\foreach \y in {1.3,1.65,2.,...,3.4}{
			\draw[fill=white] (-0.02+\x,\y) circle [radius=0.11, color=white];
		}
	}
			
	\end{tikzpicture}
	\caption{The pictures illustrate $A_{\text{ref}}(x)$ (left) taking two values in the computational domain, random defects (center), and domain mapping of the reference (right).}
	\label{fig1}
\end{center}
\end{figure}

In this paper, we study elliptic problems with diffusion coefficients that are perturbations of a single reference diffusion coefficient. We consider the following Dirichlet type problem, which we will refer to as the strong form of the (inhomogeneous) \emph{perturbed problem}: find $\bar{u}$ such that
\begin{equation}
\begin{aligned}
  \label{eq:reference_problem_classic}
	- \nabla \cdot A  \nabla \bar{u} &= f, \qquad \text{ in }\Omega, \\
	 \bar{u} & = g, \qquad  \text{ on }\Gamma,
%	 n \cdot A  \nabla u &= 0, \qquad  \text{ on }\Gamma_{\text{N}},
\end{aligned}
\end{equation}
on a bounded polygonal/polyhedral domain $\Omega \subset \mathbb{R}^d$, $d=2,3$, with boundary $\Gamma$. We assume the right hand side $f \in L^2(\Omega)$, diffusion coefficient $A \in L^{\infty}(\Omega,\mathbb{R}^{d \times d})$ is symmetric positive definite and rapidly varying, and the trace of the function $g\in H^1(\Omega)$ defines the Dirichlet boundary conditions. %The boundary $\partial \Omega$ is partitioned in disjoint subsets $\Gamma_{\text{N}}$ and $\Gamma_{\text{D}}$ for points where Neumann and Dirichlet boundary conditions apply, respectively. Further, $n$ denotes the outward normal of the boundary and $g \in H^{1}(\Omega)$ with a trace on $\Gamma_{\text{D}}$ which defines the Dirichlet boundary values.
We further consider $A$ and $f$ to be a perturbations of a reference diffusion coefficient and a reference right hand side, respectively,
\begin{equation}
A(x)\approx A_{\text{ref}}(x), \qquad f(x)\approx f_{\text{ref}}(x).
\end{equation}
The aim of the paper is to reuse computations made with the reference quantities $A_\text{ref}$ and $f_\text{ref}$ in a reliable way when solving a set of perturbed problems. The perturbations we study are local defects and domain mapppings. Figure \ref{fig1}  illustrates the reference coefficient (left), random defects (center), and domain mapping to the physical diffusion coefficient (right).

Perturbations of the diffusion coefficient in elliptic problems have been studied extensively. This work was inspired by a series of papers by Le Bris and coworkers on weakly random homogenization \cite{leb,le2014}, where they consider weakly random coefficient problems, similar to the once illustrated in Figure \ref{fig1}, using the multiscale finite element method \cite{MsFEM} for the spatial discretization. This allows the authors to consider both rapidly varying and perturbed diffusion coefficients. There are several works on partial differential equations posed on random domains using domain mapping including \cite{Harbrecht,castrillon}. In the context of perturbed coefficients, which is considered in this paper, domain mapping is instead used to transform a reference coefficient to the perturbed  coefficient. A perturbation in position of the material distribution is transformed back to a difference in value of the reference diffusion coefficient, which is advantageous from a numerical perspective. 

Multiscale methods have been a vibrant area of research for decades \cite{MsFEM,EEnq,Gamblets,LiptonBabuska}. The main idea is to solve local fine scale problems to compute an improved basis which is used to solve a global coarse-scale problem. These techniques are often parallel by construction. One method that has proven to give accurate results also for non-periodic diffusion coefficient is the Localized Orthogonal Decomposition (LOD) method \cite{MP14}. LOD is based on an orthogonal split of the solution space, using the scalar product induced by the weak form of the problem \eqref{eq:reference_problem_classic}. In recent years it has been improved and reformulated. In \cite{elf} a Petrov--Galerkin version of the method was presented and analyzed. PG-LOD has the advantage that the assembly of the modified stiffness matrix is much faster than for the original method. Concerning the implementation of the LOD, a detailed algebraic overview has been given in \cite{EPMH16}. In the recent work
\cite{hell2017} a sequence of problems with similar coefficients are considered, with applications in time dependent diffusion problems. This approach is also useful for studying perturbations of the type presented in Figure \ref{fig1}, see \cite{Tim}.

In this paper we apply the PG-LOD methodology introduced in \cite{hell2017} to solve elliptic problems with perturbed diffusion coefficient. The PG-LOD method allows for local recomputation of basis functions to handle perturbations in the data from defects and domain mappings. We derive error indicators to decide where recomputations are necessary. In this way we can efficiently simulate a vast number of perturbations by mainly solving (upscaled) coarse-scale problems.

The paper is organized as follows. In Section 2 we formulate the problem and the types of perturbations that are considered in this paper. In Section 3 we present the proposed numerical method based on PG-LOD. In Section 4 we derive error bounds and in Section~\ref{sec:implementation} we discuss implementation, memory consumption and parallelization. Finally in Section 6 we present numerical examples.

\section{Problem formulation}
\label{sec:problem_formulation}

We assume that the coefficient matrix $A \in L^{\infty}(\Omega,\mathbb{R}^{d \times d})$ is symmetric and uniformly elliptic such that
\begin{align}
	0 < \alpha &:= \essinf\limits_{x \in \Omega} \inf_{v \in \mathbb{R}^d \setminus \set{0}} \frac{\left( A(x)v \right) \cdot v}{ v \cdot v}, \\
	\infty > \beta &:= \esssup\limits_{x \in \Omega} \sup_{v \in \mathbb{R}^d \setminus \set{0}} \frac{\left( A(x)v \right) \cdot v}{ v \cdot v}.
\end{align}
We further let $f\in L^2(\Omega)$ and $g \in H^{1}(\Omega)$ with a trace on $\Gamma$ which defines the Dirichlet boundary values. We introduce a function space $V$ where we seek a solution of equation \eqref{eq:reference_problem_classic} posed on variational form. In this paper we primarily consider a conforming finite element space
$$V:=V_h\subset H^1_0(\Omega)=\{v\in H^1(\Omega)\ |\  \text{tr}(v)=0\},$$
defined on a computational mesh $\mathcal{T}_h$ that is assumed to be fine enough to resolve the variations in the diffusion coefficients well. However, we may as well choose $V=H^1_0(\Omega)$ and the analysis presented in the paper will still go through. To simplify the notation we stick to the notation $u\in V=V_h$. On weak form we get: find $u \in V$ such that
\begin{equation}
  \label{eq:perturbed_problem}
  a(u,v) =  F(v) - a(g, v)
\end{equation}
for all $v \in V$, where
\begin{equation}
 a(u,v) := \int_{\Omega}^{} \left(  A  \nabla u \right) \cdot  \nabla v, \qquad {F}(v) := \int_\Omega {f} v.
\end{equation}
The Lax--Milgram Lemma guarantees existence and uniqueness of a solution $u \in V$. The full solution including the boundary data is given by $u+g$. As mentioned above, \eqref{eq:perturbed_problem} will be referred to as the \emph{perturbed problem}.

To motivate why we cannot simply replace the perturbed coefficient with the reference coefficient, we formulate an artificial problem based on the reference coefficient and right hand side: find $u_{\text{ref}}\in V$, such that for all $v \in V$,
\begin{equation}
  \label{eq:reference_problem}
 \int_\Omega A_\text{ref}\nabla (u_\text{ref} + g)\cdot\nabla v\,dx =  \int_\Omega f_{\text{ref}} \, v.
\end{equation}

The error between $u_{\text{ref}}$ and $u$ can be bounded in the energy norm $\vertiii{\cdot} :=  a(\cdot,\cdot)^{1/2} = \norm{  A^{1/2}  \nabla \cdot}_{L^2(\Omega)}$ in the following way, 

\begin{equation}
  \label{eq:apriori_error}
  \begin{split}
  \vertiii{ u - u_\text{ref}}^2 &\leq ( A\nabla  u- A_\text{ref}\nabla u_\text{ref},\nabla (u-u_\text{ref}))+( A_\text{ref}\nabla  u_\text{ref}-  A\nabla u_\text{ref},\nabla ( u- u_\text{ref})) \\
  &= (f - f_{\text{ref}},  u- u_\text{ref}) + ((A_\text{ref}- A)\nabla (u_\text{ref} + g),\nabla(u-u_\text{ref}))\\
  &\leq \left( \frac{C_{\text{P}}}{\alpha^{1/2}} \norm{f - f_{\text{ref}}}_{L^2(\Omega)} + \frac{C_{\text{P}}}{\alpha^{3/2}}\| A_\text{ref}- A\|_{L^\infty(\Omega)}\|f_{\text{ref}}\|_{L^2(\Omega)} \right) \vertiii{ u - u_\text{ref}},
  \end{split}
\end{equation}
where $C_\text{P}$ is the Poincar\'{e} constant for $\Omega$. This error bound suggests that even local perturbations in the structure
of the coefficient or the right hand side, e.g.~by a defect or shift, may lead to very poor accuracy.
This occurs for example if we consider a problem with a highly conductive thin channel in the diffusion coefficient and a right hand side $f$ which has support inside the channel. If the channel is moved slightly so that the support of the right hand side is now outside the channel the solution will be have very differently. The error with respect to perturbations in $f$ is less severe since it is measured in the $L^2$-norm.

For this reason it is not clear how to reuse computations for the standard finite element method in a reliable way. In this paper we will treat this difficulty using a multiscale approach where solutions to localized subproblems, based on the reference coefficient, can be reused when solving for the perturbed problem.

\subsection{Perturbations}
\label{sec:perturbations}
To simplify the presentation we consider perturbations of coefficients that only takes the two values $1$ and $0 < \alpha < 1$. We emphasize that this is not necessary for the proposed method to work or for theory to hold. However, it highlights the application to composite materials which has inspired this work. Let $\Omega_{1}, \Omega_{\alpha} \subseteq \Omega$ be two disjoint subdomains of $\Omega$ with $\Omega_{1} \cup \Omega_{\alpha} = \Omega$. Let $A_{\text{ref}}$ be defined by
\[
	A_{\text{ref}} = \rchi_{\Omega_{1}} + \alpha \rchi_{\Omega_{\alpha}},  
\] 
where $\chi$ is an indicator function.

\begin{figure}
\begin{center}
	\begin{tikzpicture}[scale=0.6, every node/.style={transform shape}, node distance=1cm]
	 %nodes
	 \node[] (market) {};
	 % We make a dummy figure to make everything look nice.
	 \node[above=of market] (dummy) {};
	 \node[punkt, right=of dummy] (t) {};
	 \node[punkt, left=of dummy] (g) {}
	   edge[pil,->, bend left=0] node[auto] {$y = \psi(x)$} (t);
	  \node [left= .5cm of g]{$A_\text{ref}(x)$};
	  \node [right= .5cm of t]{$A_\text{ref} \circ \psi^{-1}(y)$ };
 	 \foreach \x in {-11.45,-11.1,-10.75,...,-7.6}{
 	 \foreach \y in {-0.45,-0.1,0.25,...,3.4}{
 	 		\draw[fill=white] (\x+6.45,\y-0.05) circle [radius=0.11, color=white];
 	 	}
 	}
	\foreach \x in {1.55,1.9,2.25,...,5.4}{
	\foreach \y in {0.1,0.45,...,3.8}{
			\draw[fill=white] (\x-0.05+4* 1.7 * \y / 3.85 *\x/3.85-4* 1.7 * \y / 3.85 * 1.55/3.85 -4* 1.7 * \y / 3.85 * \x^2/3.85^2 +4* 1.7 * \y / 3.85 * 2* 1.55 * \x / 3.85^2 -4* 1.7 * \y / 3.85 * 1.55^2/3.85^2 - 4* 1.7 * \y^2 / 3.85^2 *\x/3.85 + 4* 1.7 * \y^2 / 3.85^2 *1.55/3.85 + 4* 1.7 * \y^2 / 3.85^2 *\x^2/3.85^2 -4* 1.7 * \y^2 / 3.85^2 * 2* 1.55 * \x / 3.85^2 + 4* 1.7 * \y^2 / 3.85^2 *1.55^2/3.85^2 ,\y - 0.6)  circle [radius=0.11, color=white];
		}
	}
	\end{tikzpicture}
	\caption{Illustration of a domain mapping $\psi$ with $D \equiv 0$.}
	\label{neu}
\end{center}
\end{figure}
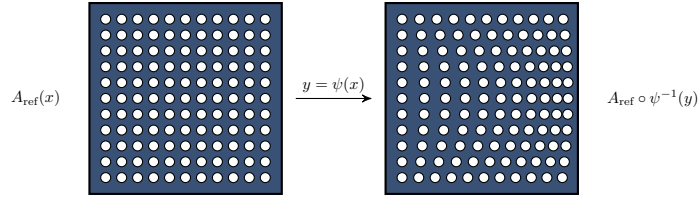

We formalize the two types of perturbations that we consider (see Figure \ref{fig1}) by introducing a defect perturbation $D$ and a domain mapping $\psi$.
A perturbation from a defect can be expressed by $D = (1 - \alpha) \rchi_{\omega}$ where $\omega \subseteq \Omega_{1}$ and that $A_{\text{ref}} - D$ can be considered as the perturbed coefficient.
For shift perturbations, we assume that the domain mapping perturbation can be described as a variable
transformation with a perturbation function
$\psi:\Omega\rightarrow\Omega$ which maps the reference coefficient
(expressed in $x$-coordinates) to a mapped coefficient (expressed
in $y$-coordinates). We assume that $\psi$ maps the boundary to itself (i.e.\ $\Gamma = \{\psi(x)\,:\,x\in\Gamma\}$) and that it is a one-to-one
mapping in $\Omega$. Figure \ref{neu} provides an example of a variable
transformation. We denote the corresponding Jacobi matrix
\begin{equation*}
  \mathcal{J}_{ij}(x)=\left[\frac{\partial \psi_i}{\partial x_j}\right],
\end{equation*}
for $1\leq i,j\leq d$, and assume it to be bounded with bounded
inverse for a.e.\ $x\in \Omega$.  The two perturbation types can
either be combined, or be considered individually by letting
$D \equiv 0$ or $\psi = \id$.

Using $A_{\text{ref}}$, $D$ and $\psi$, we formulate the \emph{mapped problem} in the $y$-variable with $y = \psi(x)$ by
\begin{equation} \label{eq:domain_mapping_problem}
  \begin{split}
    - \nabla_y \cdot A_y \nabla_y u_y & = f_y, \qquad \text{ in }\Omega, \\
    u_y & = g_y, \qquad \text{ on } \Gamma,
\end{split}
\end{equation}
where the coefficient is defined by
$A_y = (A_\text{ref}-D) \circ \psi^{-1}$, and the derivatives have
been distorted accordingly. The $y$-variable corresponds to the
physical spatial variable in a typical situation.  Depending on the
physics being modeled, either the mapped right hand side
$f_y \in L^2(\Omega)$ or the perturbed $f$ (below) can be considered
given. It makes no difference for the development of the numerical
method, but may affect the choice of $f_{\text{ref}}$. The Dirichlet
boundary value $g_y \in H^{1/2}(\Gamma)$ is defined only on the
boundary $\Gamma$, which is mapped to
itself. %Thus, the domain mapping has no influence on the Dirichlet boundary conditions.
The solution in the perturbed domain is denoted $u_y(y)$.

%\begin{equation} \label{eq:domain_mapping_problem}
%	\begin{split}
%		- \nabla_y \cdot \left( A_\text{ref}(\psi^{-1}(y))  \nabla_y \bar{u}(y) \right) &= f_\text{ref}(\psi^{-1}(y)), \qquad \text{ in }\Omega, \\
%		 \bar u(y) &= g(y), \qquad \text{ on } \Gamma_{\text{D}}\\
%n\cdot A_\text{ref}(\psi^{-1}(y))\nabla \bar{u}(y)=0, \qquad \text{ on } \Gamma_{\text{N}}.
%	\end{split}
%\end{equation}

Next, we use the mapped problem to define the perturbed problem in equation \eqref{eq:perturbed_problem}. The gradient operator $\nabla_y$ in the mapped domain can be expressed in $\nabla_x = \nabla$ by
\begin{equation}
  \nabla_y v(x) := \left[ \diff{v}{y_i}(x)  \right]_i = \left[ \sum_{j}^{}\diff{v}{x_j}(x) \diff{x_j}{y_i}(x) \right]_i = \mathcal{J}^{-T}(x)  \nabla_x v(x).
\end{equation}
%Further, the mapped outward boundary normal is $n_y = \mathcal{J}^{-T} n$.
Based on the elliptic operator in equation \eqref{eq:domain_mapping_problem} we define the perturbed bilinear form for the mapped problem as %use the weak formulation of \eqref{eq:domain_mapping_problem} to obtain
\begin{equation*}
  \begin{split}
    a(v,w) &= \int_{{\Omega}}^{} \left( (A_\text{ref}-D) \circ \psi^{-1} \right )  \nabla_y \left(v \circ \psi^{-1} \right) \cdot  \nabla_y \left(w \circ \psi^{-1} \right) \, \dd y\\
    &= \int_{\Omega}^{} \det(\mathcal{J}) \mathcal{J}^{-1} (A_\text{ref}-D) \mathcal{J}^{-T}  \nabla_x v \cdot  \nabla_x w \, \dd x
  \end{split}
\end{equation*}
and the corresponding linear functional
\begin{equation*}
  {F}(w) = \int_{ \Omega} f_y \left(w \circ \psi^{-1}\right) \,\dd y = \int_{\Omega} \det(\mathcal{J}) \left(f_y\circ \psi\right)  w\, \dd x.
\end{equation*}
We see that this now fits the formulation of the perturbed problem \eqref{eq:perturbed_problem} with
\begin{equation*}
  A=\det(\mathcal{J}) \mathcal{J}^{-1} (A_\text{ref}-D) \mathcal{J}^{-T}, \qquad f = \det(\mathcal{J}) \left(f_y\circ \psi\right), \qquad   g|_{\Gamma} = g_y \circ \psi,
\end{equation*}
and mapped solution $u_y = (u + g) \circ \psi^{-1}$. For the problem to be well posed we assume $\mathcal{J}$ and $\mathcal{J}^{-1}$ to be bounded almost everywhere and $A$ to be symmetric positive definite. We note that the perturbed coefficient $A$ can be computed from the reference coefficient by means of the Jacobian matrix and that the Dirichlet boundary value function $g \in H^1(\Omega)$ can be chosen arbitrarily in the interior of $\Omega$.
We emphasize that the domain mapping transforms a shift defect to a change-in-value perturbation. For many coefficients this is advantageous as seen in equation \eqref{eq:apriori_error}, where now the $L^{\infty}$-norm can be expressed entirely in terms of how much $\mathcal{J}$ differs from identity. The domain mapping covers continuous (possibly global) perturbations while defects cover discontinuous (often local) perturbations. 	With respect to Figure \ref{fig1} we see that the middle picture corresponds to $\psi = \id$ and the right picture to $D \equiv 0$. 

%\todo{We should define exactly what conditions we need on $\psi$ for this to hold}
%Finally, we emphasize that the domain mapping, which is a coefficient shift perturbation, has been replaced by a change-in-value perturbation by mapping the domain back to the reference domain, see Figure \ref{fig2.2}. For many coefficients this is advantageous as seen in equation \eqref{eq:apriori_error}, where now the $L^{\infty}$-norm can be expressed entirely in terms of $\mathcal{J}$.

%In the next few sections we develop a method that reuse corrector computations performed with the reference coefficients on problems with a perturbed coefficients. We will focus on the general perturbed problem \eqref{eq:perturbed_problem} expressed in $A$, $u$, $f$, and $g$, where the perturbed coefficient stems from either a domain mapping or defects.

\begin{bemerkung}[Discretization of $\psi$]
In the numerical examples in Section 6 we let $V=V_h$ be the space of quadrilateral finite elements and $\psi$ to be a linear combination of the corresponding bilinear shape functions. This leads to an isoparametric finite element formulation. Isoparametric finite elements are for instance used to guarantee accuracy when solving problems on curved domains. Here the aim is instead to map the perturbed diffusion coefficient to a reference. The theoretical justification for using bilinear domain mappings $\psi$ follows directly from the theory of isoparametric finite elements, see \cite{Ciarlet,BrennerScott}.
\end{bemerkung}

%\fredrik{This paragraph seems a bit off here.} The aim of this work is to present an efficient method for solving a large number of perturbed problems. To avoid global recomputation we introduce multiscale techniques for which local updates are possible.

%In the applications, we intend to solve several problems of the same type, based on the same reference coefficient. Therefore, a reduction of the computational costs is crucial. In terms of the standard FEM, we either compute the stiffness matrix with respect to the reference coefficient $\hat A$ or to the perturbed coefficient $A$. According to Lemma \ref{ayay}, approximating a perturbed problem with the reference solution causes a large error, whereas a recomputation is much more expensive and does not rely on the reference configuration. This fact is not satisfying in terms of computational complexity. The next section is devoted to an algorithm that resolves this issue.

\section{Adaptive PG-LOD method}
In this section, we develop the proposed method that reuses reference computations peformed using the reference coefficient $A_{\text{ref}}$ and right hand side $f_{\text{ref}}$ in an adaptive manner to solve the perturbed problem \eqref{eq:perturbed_problem} expressed in $A$, $u$, $f$, and $g$. The perturbed data can stem from a domain mapping, defects or both. 
We consider the Petrov--Galerkin version of the LOD method as presented in \cite{elf} and briefly derive it in Sections~\ref{sec:prel}--\ref{sub:localization}. Section~\ref{sec:error_indicators} defines the computable error indicators and Section~\ref{sec:adaptive} presents the method that reuses reference computations while balancing the error. Recently, a similar approach was applied to problems with time dependent diffusion, see \cite{hell2017}.

\subsection{Preliminaries}\label{sec:prel}
Let $\TT _H$ be a coarse, shape regular, conforming mesh family of the domain $\Omega$. We denote the maximum diameter of an element in $\TT _H$ with $H$ and $\mathcal{N}$ the set of all corresponding interior nodes of the mesh $\TT _H$.
We let
\[
	V_H := V \cap \mathcal{P}_1(\TT _H),
\]
where $\mathcal{P}_1(\TT _H)$ denotes the space of $\TT _H$-piecewise affine functions that are continuous on the domain $\Omega$. Since our full space $V$ is a conforming finite element space we assume that the meshes and spaces are nested $V_H\subset V$. However, it is possible, with a minor modification of the proposed method, to violate this condition and still get convergence, see \cite{MaPe15}.

\begin{figure}[h]
	\centering
	\includegraphics[scale=0.3]{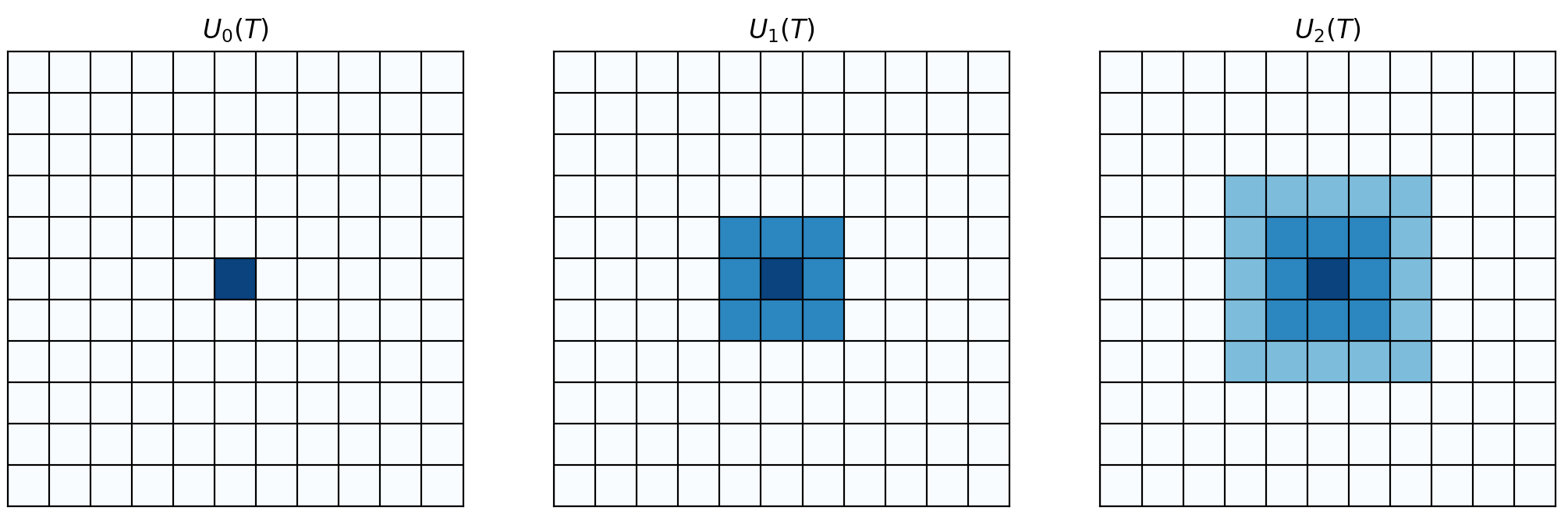}
	\caption{Patches for a coarse mesh element $T \in \TT _H$.}
	\label{patches}
\end{figure}

We introduce the concept of element patches as they will be used in the definition of the interpolation and for the localization of the PG-LOD method. For arbitrary $\omega \subseteq \Omega$ and $0 \leq k \in \mathbb{N}$, we define coarse grid patches $U_k(\omega) \subset \Omega$ by
\begin{align*}
	U_0(\omega) &= \omega,  \\
	U_{k+1}(\omega) &= \bigcup \set[T \in \TT _H]{ \overline{U_k(\omega)} \cap \overline{T} \neq 0}.
\end{align*}
If $\omega = \set{x}$, for a node $x \in \mathcal{N}$, we call $U_k(x)$ a $k$-layer nodal patch. For $\omega = T$, where $T \in \TT _H$, we call $U_k(T)$ a $k$-layer element patch, see Figure \ref{patches}

\subsection{Multiscale decomposition}
\label{sub:multiscale_splitting}
Fine scale features that occur in the solutions are not captured in the space $V_H$. We characterize the fine scale parts of $V$ as the kernel of an linear surjective (quasi-)interpolation operator $\mathcal{I}_H : V \to V_H$ that maps a function $v \in V$ to a function $v_H \in V_H$ in the coarse FE space. Let $\TT_x = \{T \in \TT_H\,:\, x \in \overline{T}\}$ be the set of elements neigboring $x$. We let the interpolation operator $\mathcal{I}_H:V \to V_H$ be an $L^2$-projection to the broken finite element space in composition with averaging in the nodes. In other words,
\[
  \mathcal{I}_Hv := \sum_{x \in \mathcal{N}} \zeta_x \lambda_x,
\]
where coefficients $\zeta_x$ are determined as follows. Let the operator $P_T$ be the element $L^2$-projection $P_Tv \in V_H  \big|_{T}$ such that
\[
  \int_{T}(P_Tv)w_H = \int_{T}vw_H,
\]
for all $w_H \in V_H \big|_{T}$, then
\[
  \zeta_x = \operatorname{card}(\TT_x)^{-1} \sum_{T' \in \TT_x} (P_{T'}v)(x).
\]
This interpolation operator is linear, continuous and its restriction to $V_H$ is an isomorphism. Furthermore it fulfills the stability result
\begin{equation}
H^{-1}_T \norm{v - \mathcal{I}_Hv}_{L^2(T)} + \norm{ \nabla \mathcal{I}_Hv}_{L^2(T)} \leq C_{\mathcal{I}_H} \norm{ \nabla v}_{L^2(U_1(T))}, \label{Ia}
\end{equation}
for every $v \in V$ and $T \in \TT _H$, with a generic constant $C_{\mathcal{I}_H} >0$, see \cite{Pe15}. We refer to \cite{MH16, PS} for other possible choices of interpolation operators.

We let the kernel of $\mathcal{I}_H$, \[
	V^{\text{f}} = \ker(\mathcal{I}_H) = \set[v \in V]{\mathcal{I}_H(v)=0}
\]
define the fine scales of the space $V$. Since $\mathcal{I}_H$ is a projection, this allows for the split $V = V_H \oplus V^{\text{f}}$.
%This procedure allows for covering all features of $V$ that are no longer contained in $V_H$ and we yield the multiscale splitting
We define a correction operator $\QQ v \in V^{\text{f}}$, for a given $v \in V$, to be the solution of
\[
  a(\QQ v,v^{\text{f}}) = a(v,v^{\text{f}}),
\]
for all $v^{\text{f}} \in V^{\text{f}}$,
%In particular, for the subset $V_H \subset V$, $\QQ $ denotes a fine scale projection to $V^{\text{f}}$.
and define the multiscale space $V^{\text{ms}} := V_H - \QQ V_H$.

For any $v^{\text{f}} \in  V^{\text{f}}$ and $v^{\text{ms}} \in V^{\text{ms}}$, we observe $a(v^{\text{ms}},v^{\text{f}}) =0$. This leads to the orthogonal decomposition with respect to the $a$-scalar product, $V = V^{\text{ms}} \oplus_a V^{\text{f}}$.
 Right hand correction can be used to improve accuracy in LOD based methods, see e.g.~\cite{malhen2014, hell2017}. We define this correction by $\RR f \in V^{\text{f}}$ such that, for all $v^{\text{f}} \in V^{\text{f}}$,
 \begin{align*}
 	a(\RR f,v^{\text{f}}) = \int_\Omega fv^{\text{f}}.
 \end{align*}
We now derive the Petrov--Galerkin LOD method for the perturbed problem, see also \cite{elf,GallPet}. We use $V = V_H \oplus V^{\text{f}}$ to decompose \eqref{eq:perturbed_problem} into two equations 
\begin{align}
  a(u_H + u^{\text{f}},v_H) &= F(v_H) - a(g, v_H), \label{derive:LOD_1} \\
  a(u^{\text{f}},v^{\text{f}}) &= F(v^{\text{f}}) - a(g, v^{\text{f}}) - 	a(u_H,v^{\text{f}}), \label{derive:LOD_2}
\end{align}
for all $v_H \in V_H$ and $v^{\text{f}} \in V^{\text{f}}$.

With the definitions of $\QQ$ and $\RR$ we obtain from \eqref{derive:LOD_2} that $u^{\text{f}} = -\QQ u_H + \RR f - \QQ g$. Plugging this into \eqref{derive:LOD_1} gives
\begin{align}
a(u^{\text{ms}} ,v_H) =  F(v_H) - a(\RR f, v_H) - a(g - \QQ g, v_H), \label{derive:LOD_3}
\end{align}
for all $v_H \in V_H$. Hence, solving \eqref{derive:LOD_3} gives the exact solution $u = u^{\text{ms}}+ \RR f - \QQ g$.

\begin{figure}
	\begin{subfigure}[b]{0.33\textwidth}
		\includegraphics[width=\textwidth]{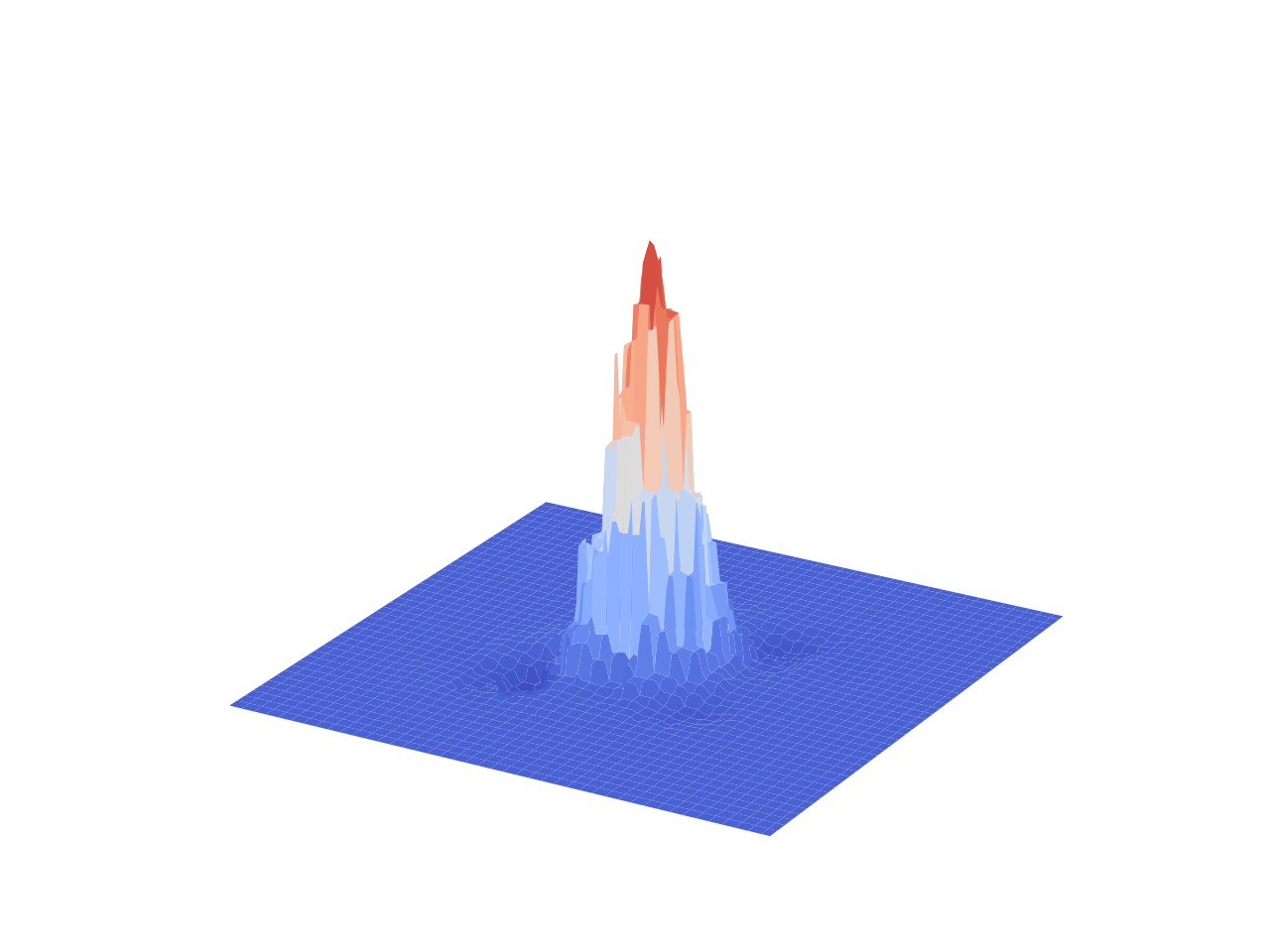}
		\caption{$\lambda_x - \QQ \lambda_x$.}
	\end{subfigure}
	\begin{subfigure}[b]{0.33\textwidth}
		\includegraphics[width=\textwidth]{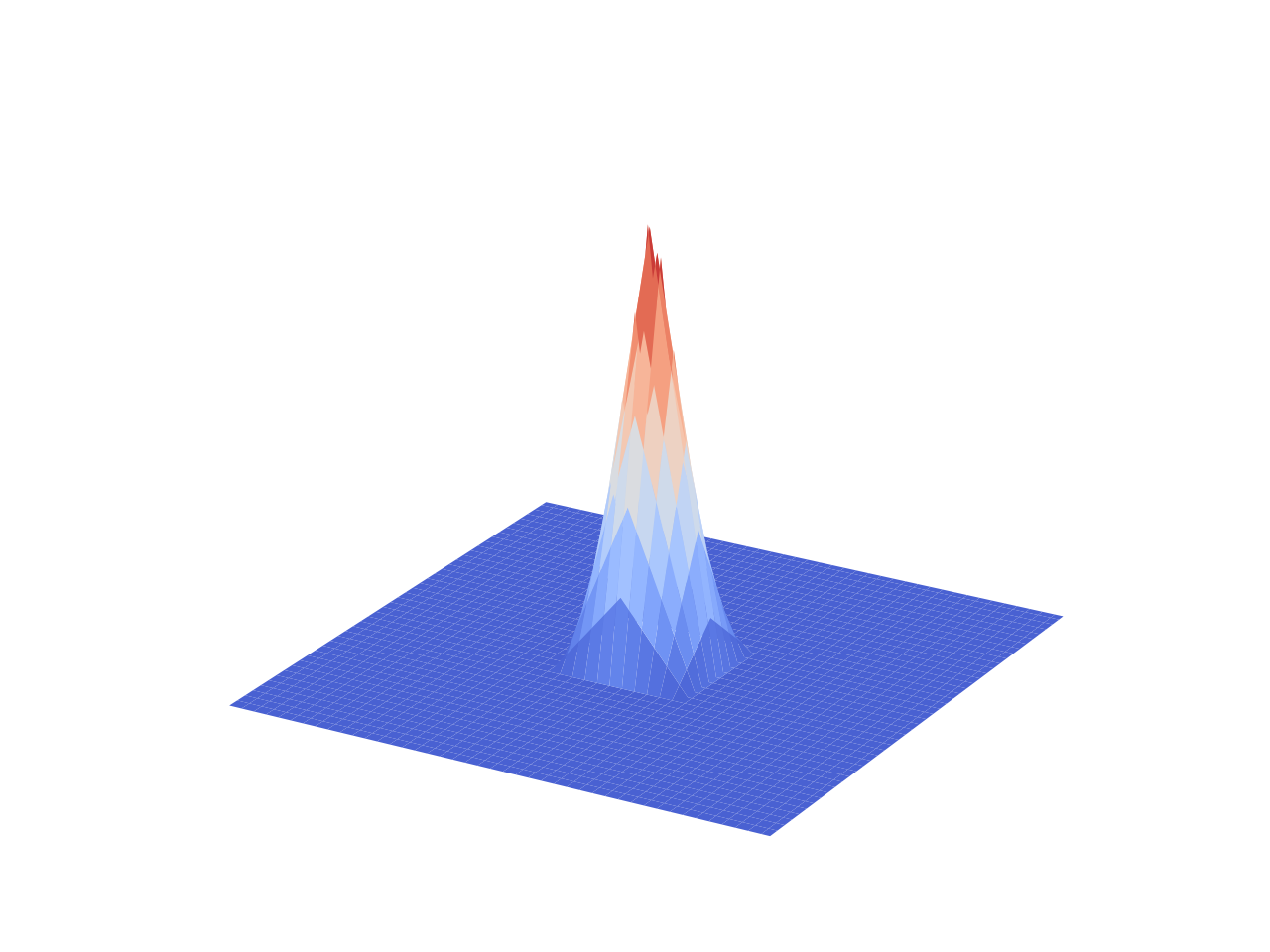}
		\caption{$\lambda_x$.}
	\end{subfigure}
	\begin{subfigure}[b]{0.33\textwidth}
		\includegraphics[width=\textwidth]{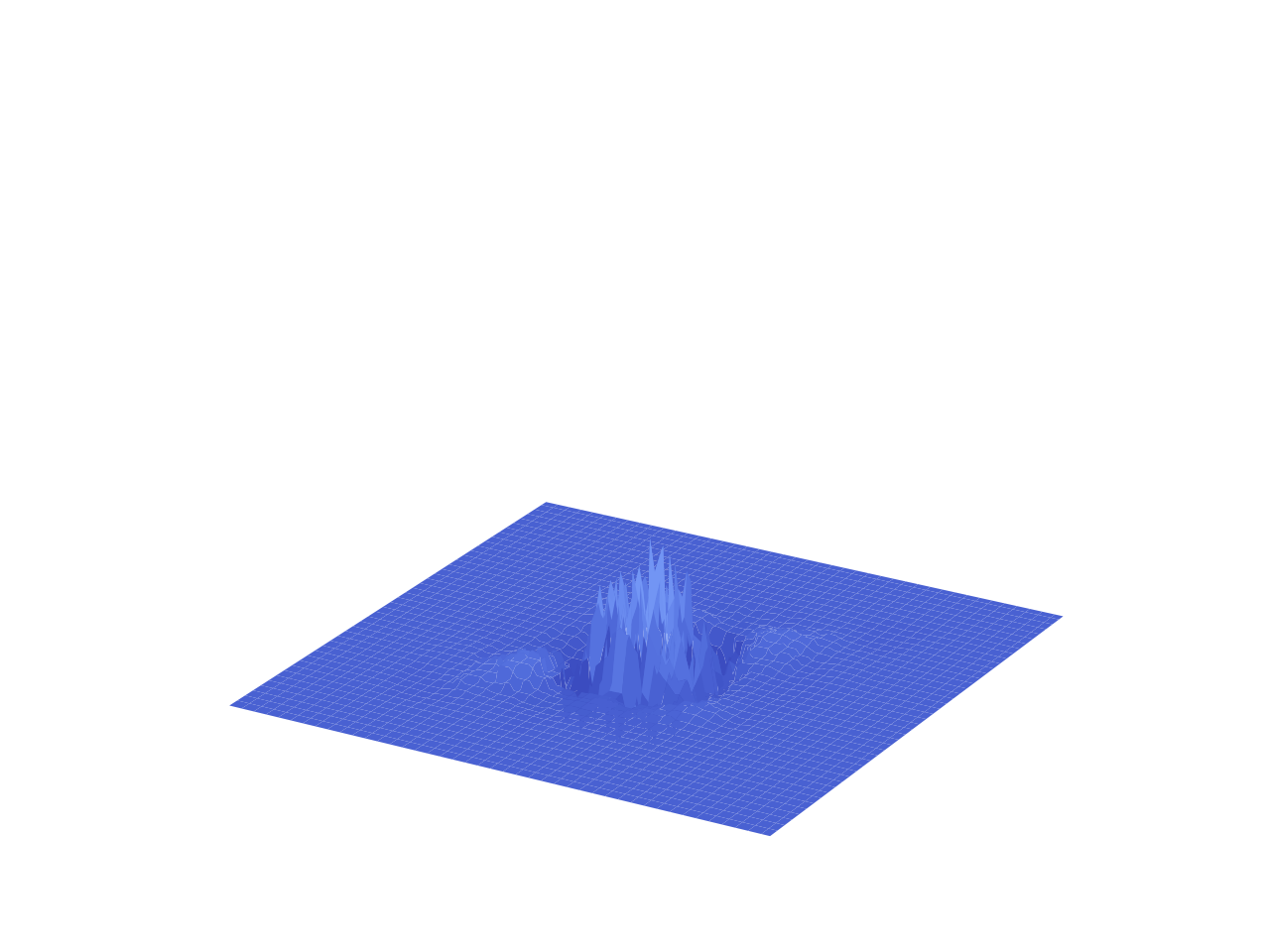}
		\caption{$\QQ \lambda_x$.}
	\end{subfigure}
	\caption{Basis function of $V_{H}^\text{ms}$ and its decomposition for $x \in \mathcal{N}$.}
	\label{correctorplots}
\end{figure}

In order to use the multiscale space in a practical implementation, we need a computable basis.
Since $V^{\text{ms}}$ and $V_H$ have equal dimensions, it suffices to apply the fine scale corrector $\QQ $ on every single basis function $\lambda_x$ of $V_H$ to obtain a corrected basis, i.e.
% and incorporate it to the underlying basis function. Therefore, for every node $x \in \mathcal{N}$, we define $\phi_x := \QQ \lambda_x \in V^{\text{f}}$, which means that $\phi_x \in V^{\text{f}}$ solves, for every $w \in V^{\text{f}}$
%\begin{equation}
%	a(\phi_x,w) = a(\lambda_x,w). \label{fix}
%\end{equation}
\[
	\set[\lambda_x - \QQ \lambda_x]{x \in \mathcal{N}}.
\]
We note that a global fine scale computation for each node is necessary to compute all $\QQ \lambda_x$ which is computationally expensive. However, Figure~\ref{correctorplots} suggests that the computation can be localized to a small area around the support of the original basis function. It was shown in \cite{MP14} that the corrected basis functions decay exponentially, and that localized computations are possible.

%Due to that, Målqvist and Peterseim investigated the behavior of these globally supported correctors. A view of Figure \ref{correcorplots} reveals that the correctors indeed have an exponential decay outside an area of their node. They presented a proof in \cite{MP14} and used it as the justification for truncating the corrector operator. Thus, it is computed on only a patch around the affected node $x \in \mathcal{N}$ and we gain a method that is localized on every coarse grid patch.
\begin{bemerkung}[Mixed boundary conditions]
It is also possible to have mixed Neumann and Dirichlet boundary conditions in the PG-LOD method, see \cite{malhen2014}. With mixed boundary conditions we also need to compute fine scale corrections for the inhomogeneous Neumann boundary data (as we do for the Dirichlet data) to get optimal convergence rate.
\end{bemerkung}

\subsection{Localized multiscale method}
\label{sub:localization}
The fine scale space $V^{\text{f}}$ can be restricted to patches $U_k(\omega)$ with the intuitive definition, for $\omega \subseteq \Omega$ and $k \in \mathbb{N}$,
\[
	V^{\text{f}}(U_k(\omega)) := \set[v \in V^{\text{f}}]{ v  \big|_{\Omega \setminus U_k(\omega)}^{} = 0}.
\]
These local fine scale patches enable the truncation of the corrector. We define localized element correction operators $\QQ _k$ and $\RR_k$ by
\[
	\QQ _k v := \sum_{T \in \TT _H}^{}\QQ _{k,T} v, \qquad \text{and} \qquad \RR_kf := \sum_{T \in \mathcal{T}_H}^{} \RR_{k,T}f,
\]
where $\RR_{k,T},\QQ _{k,T}\,:\, V \to V^{\text{f}}(U_k(T))$ solves
\begin{align*}
	a(\QQ _{k,T} v,v^{\text{f}}) &= \int_T A\nabla v\cdot\nabla v^{\text{f}}, \\
	a(\RR_{k,T} f,v^{\text{f}}) &=  \int_T f \, v^{\text{f}},
\end{align*}
for all $v^{\text{f}} \in V^{\text{f}}(U_k(T))$.
We construct a localized multiscale space $V_{k}^{\text{ms}}$ using $V_H$ and the local correctors
\[
	V_{k}^{\text{ms}} := V_H - \QQ _kV_H.
\]
This space is spanned by $\set{\lambda_x - \QQ _k \lambda_x}_{x \in \mathcal{N}}$.

We formulate the localized version of the Petrov--Galerkin LOD in \eqref{derive:LOD_3}: find $u^{\text{ms}}_k \in V_{k}^{\text{ms}}$ such that for all $v \in V_H$,
\begin{align}
  \label{eq:PGLODref}
  a(u^{\text{ms}}_k,v) = F(v) - a(\RR_kf, v) - a(g-\QQ _kg, v)
\end{align}	
where the full approximation of $u$ is
\begin{equation}\label{eq:localized_uk}
  u_k = u^{\text{ms}}_k + \mathcal{R}_kf - \QQ_k g = u^{H}_k - \QQ_k u^H_k + \mathcal{R}_kf - \QQ_k g.
\end{equation}
 The main reason for using a Petrov--Galerkin formulation is that it avoids the expensive computation of products between corrected basis functions without losing convergence order \cite{elf}. 
% We make the choice to omit the right hand side correction for $f$, but include it for $g$. The error introduced from omitting the correction for $f \in L^{2}(\Omega)$ is proportional to $H$, while leaving the $H^{-1}(\Omega)$ functional $a(g,\cdot)$ uncorrected would yield of order one. See e.g.\ \cite{MH14} for a detailed discussion on this topic.

 The exponential decay of the correctors yield the following error bounds (see \cite{MP14}) of the localized correctors in terms of $k$:
 \begin{equation}\label{eq:decay}
     \vertiii{(\QQ - \QQ_k) v} \lesssim k^{d/2} \theta^k \vertiii{v}, \qquad\text{and}\qquad
     \vertiii{(\RR - \RR_k) f} \lesssim k^{d/2} \theta^k \|f\|_{L^2(\Omega)},
 \end{equation}
 where the notation $a\lesssim b$ means $a\leq C b$ with a constant $C$ independent of $H$, $k$, and TOL (which is used in the coming sections). The well-posedness of \eqref{eq:PGLODref} was studied previously in \cite{elf, hell2017} and appears to be conditioned on sufficiently large $k$ in general. (This condition will be revisited in the proof of Theorem \ref{theorem}).
 Furthermore, from \cite[Section 4.2]{hell2017} we obtain an error bound for $u - u_k$ which reads  
 \begin{equation} \label{error_bound_1}
   \begin{split}
     \vertiii{u-u_k} & \lesssim \vertiii{(\QQ - \QQ_k)(\IH u + g)} + \vertiii{(\RR - \RR_k)f} \lesssim  k^{d/2}\, \theta^k \left(\|f\|_{L^2(\Omega)} + \vertiii{g} \right).
   \end{split}
\end{equation} 

\subsection{Error indicators}
\label{sec:error_indicators}
Since the aim of this paper is to reuse local corrector computations performed with the reference coefficients $A_{\text{ref}}$ and $f_{\text{ref}}$,
we need a notation for a modified $\QQ _{k,T}$ 
and $\RR _{k,T}$, computed using $A_\text{ref}$ instead of $A$ and $f_{\text{ref}}$ instead of $f$. We let $\QQ ^{\text{ref}}_{k,T}\,:\, V \to V^{\text{f}}(U_k(T))$ 
and $\RR ^{\text{ref}}_{k,T}\,:\, L^2(\Omega) \to V^{\text{f}}(U_k(T))$
be the solutions to
\begin{align*}
  (A_{\text{ref}}\nabla\QQ ^{\text{ref}}_{k,T} v,\nabla v^{\text{f}}) &= \int_T A_\text{ref}\nabla v \cdot \nabla v^{\text{f}}, \\ 
  (A_{\text{ref}}\nabla\RR ^{\text{ref}}_{k,T} f_{\text{ref}},\nabla v^{\text{f}}) &= \int_T f_{\text{ref}} v^{\text{f}},
\end{align*}
for all $v^{\text{f}} \in V^{\text{f}}(U_k(T))$.
In order to decide for which $T\in \TT_H$ we need to recompute the correctors ($\QQ_{k,T}$ 
and $\RR_{k,T}$) and for which we can still use the reference correctors ($\QQ ^{\text{ref}}_{k,T}$ 
and $\RR ^{\text{ref}}_{k,T}$)
we need computable error indicators. For the indicators to be useful and efficient they have to be computable and independent of $\QQ_{k,T}$ and $\RR_{k,T}$.

For every $T \in \TT _H$ we define error indicators only depending on the reference corrector and coefficients $A$ and $A_\text{ref}$. The definition of the error indicators is motivated by Lemma~\ref{lemma:indicator} and by Theorem~\ref{theorem} in Section~\ref{cha:error_analysis_for} below. See \cite{hell2017} for a similar construction.
\begin{definition}[Error indicators] \label{def:indicators}
For each $T \in \TT _H$, we define
	\begin{equation}  \label{indidef}
		\begin{split}
			E^2_{\QQ V_H,T} &:= \|A_\text{ref}A^{-1}\|_{L^\infty(T)}\sum_{\substack{T' \in \TT _H \\ T' \cap U_k(T) \neq 0}}^{} \mkern-10mu \norm{\delta}^2_{L^{\infty}(T')}
			\max_{\substack{ w |_{T}, \,w \in V_H %,\\ \|A_\text{ref}^{1/2}\nabla w\|_{L^2(T)}=1
}
} \frac{\norm{ A_{\text{ref}}^{1/2} (\rchi_T  \nabla w -  \nabla {\QQ ^{\text{ref}}_{k,T}}w)  }^2_{L^2(T')}}{\|A_\text{ref}^{1/2}\nabla w\|^2_{L^2(T)}}, \\
E^2_{f,T} &:= \|A^{-1}\|_{L^\infty(U(T))} \sup_{\substack{ v \in V}} \frac{\|(\id-\mathcal{I}_H)v\|^2_{L^2(T)}}{\|\nabla v\|^2_{L^2(U(T))}} \norm{f - f_{\text{ref}}}^2_{L^2(T)}, \\
E^2_{\RR f,T} &:= \|A_\text{ref}A^{-1}\|_{L^\infty(T)}\mkern-10mu \sum_{\substack{T' \in \TT _H \\ T' \cap U_k(T) \neq 0}}^{}\mkern-10mu \norm{\delta}^2_{L^{\infty}(T')}
\norm{ A_{\text{ref}}^{1/2} \nabla {\RR ^{\text{ref}}_{k,T}}f_{\text{ref}} }^2_{L^2(T')},  \\
			 E^2_{\QQ g,T} &:= \|A_\text{ref}A^{-1}\|_{L^\infty(T)}\sum_{\substack{T' \in \TT _H \\ T' \cap U_k(T) \neq 0}}^{}\mkern-10mu \norm{\delta}^2_{L^{\infty}(T')}
			 \norm{ A_{\text{ref}}^{1/2} (\rchi_T  \nabla g -  \nabla {\QQ ^{\text{ref}}_{k,T}}g)  }^2_{L^2(T')}, \\
		\end{split}
	\end{equation} 
	where $\rchi_T$ denotes the indicator function for an element $T \in \TT _H$ and $$\delta=A^{-1/2}(A-A_\text{ref})A_\text{ref}^{-1/2}.$$
\end{definition}

We define the square root of a symmetric positive definite matrix $ A^{1/2}$ as the unique principal square root, also positive definite. The error indicators are defined with the goal to reduce memory consumption. With the definitions above, all quantities depending on the reference coefficient and right hand side can be computed in advance. Further implementation details will be discussed in Section~\ref{sec:implementation}.

The numerical method we propose in the next subsection exploits the possibility to replace e.g.\ $\RR_{k,T}f$ with the precomputed $\RR^{\text{ref}}_{k,T}f_{\text{ref}}$, and to reuse integrals for the stiffness matrix based on $A_{\text{ref}}$ instead of using $A$. The error indicators will be used to identify which local problems that need to be recomputed and which local problems that can be computed using the reference data. 

\begin{bemerkung}[Perturbed boundary data]
We do not cover the case when $g$ is perturbed. It can
be perturbed explicitly from a reference boundary
condition (in which case a reference $g$ similar to
$f_{\text{ref}}$ needs to be introduced), or it
can be perturbed by a domain mapping $g|_{\Gamma} = g_y \circ \psi$. We remark that if
$\psi$ is the identity mapping for all points on the
boundary $\Gamma$, then $g$ (including its values in
the interior of the domain) can be picked independent
of $\psi$.
\end{bemerkung}

\begin{bemerkung}[Pure domain mapping]
In the domain mapping setting, assuming $A_{\text{ref}}$ to be scalar and $D\equiv0$, we have that
 $$\delta=(\det(\mathcal{J}) \mathcal{J}^{-1} \mathcal{J}^{-T})^{-1/2} \left( \det(\mathcal{J}) \mathcal{J}^{-1} \mathcal{J}^{-T} -\id  \right),$$ i.e.\ it is independent of $A_\text{ref}$. If the Jacobi matrix can be written as a $\epsilon$-perturbation of the identity the size of $\delta$ will be  proportional to $\epsilon$.
\end{bemerkung}

\subsection{Adaptive method}\label{sec:adaptive}
We are now ready to present the full method with adaptively updated
correctors. The main idea is to compute the perturbed correctors only
for a subset of all elements, and reuse the reference correctors for
all other elements. This means we effectively solve the problem in a
mixed multiscale space using a bilinear form that is defined as a combination
of the two coefficients.
\begin{definition}[PG-LOD with adaptively updated correctors] \label{procedure}
The proposed method follows five steps:
\begin{enumerate}
\item
  Provided reference data $(A_{\text{ref}}, f_{\text{ref}})$: Compute (for all $T\in\TT _H$) reference correctors $\QQ ^{\text{ref}}_{k,T}\lambda_x$ (for all basis functions $\lambda_x$), $\RR ^{\text{ref}}_{k,T}f_{\text{ref}}$, and $\QQ^{\text{ref}}_{k,T}g$, based on the reference coefficient $A_\text{ref}$ and reference right hand side $f_{\text{ref}}$.
\item Provided perturbed data $(A, f)$: Compute (for all $T\in\TT _H$) error indicators $E_{\QQ V_H,T}$, $E_{f,T}$, $E_{\RR f,T}$, and
  $E_{\QQ g,T}$ and mark the elements $T$ for which all
  of the following inequalities hold true,
  \begin{equation} \label{eq:local_estimates}
    \begin{aligned}
      E_{\QQ V_H,T} & \le \text{TOL}, \\
      E_{f, T} + E_{\RR f,T} & \le \text{TOL}\left(\norm{f}_{L^2(\Omega)} + \vertiii{g}\right), \\
      E_{\QQ g,T} & \le \text{TOL}\left(\norm{f}_{L^2(\Omega)} + \vertiii{g}\right). \\
    \end{aligned}
  \end{equation}
  Denote the set of marked elements by $\TT^{\text{ref}}_H\subset \TT _H$.
\item
  Compute (for all $T \in \TT _H\setminus \TT^{\text{ref}}_H$) the mixed correctors $\tilde{\QQ }_{k,T}\lambda_x$, $\tilde{\QQ }_{k,T}g$, and $\tilde{\RR }_{k,T}\tilde f$, based on the following definitions of the mixed right hand side and correcctors: $\tilde f \in L^2(\Omega)$ and
\begin{equation*}
  \tilde{f}|_T = \left\{
    \begin{array}{ll}
      f_{\text{ref}}|_T, \\
      f|_T,
    \end{array}\right.
  \quad
  \tilde{\QQ }_{k,T} = \left\{
    \begin{array}{ll}
      {\QQ }^{\text{ref}}_{k,T}, \\
      \QQ _{k,T},
    \end{array}\right. 
  \quad
  \tilde{\RR}_{k,T} = \left\{
    \begin{array}{ll}
      {\RR }^{\text{ref}}_{k,T},& \qquad \text{for }T\in \TT^{\text{ref}}_H,\\
      \RR _{k,T},& \qquad \text{otherwise}.
    \end{array}\right.
\end{equation*}
Note that only element correctors in $\TT _H\setminus \TT^{\text{ref}}_H$ need to be recomputed since the reference correctors from Step 1 will be used for the elements in $\TT^{\text{ref}}_H$.  We further let $\tilde{\QQ }_k=\sum_{T\in\TT _H}\tilde{\QQ }_{k,T}$ and $\tilde{\RR }_k=\sum_{T\in\TT _H}\tilde{\RR }_{k,T}$. 
\item Assemble the adaptively updated LOD stiffness matrix
\begin{equation} \label{eq:K_tilde}
\tilde K_{xy} =\tilde{b}(\lambda_y,\lambda_x),
\end{equation}
using the mixed unsymmetric bilinear form $\tilde b$ defined in terms of a element-wise reference $b_{T}^{\text{ref}}$ and a perturbed $b_{T}$:
\begin{align*}
  b_{T}^{\text{ref}}(v,w) & = \skal{A_\text{ref}(\rchi_T  \nabla  -  \nabla \QQ ^{\text{ref}}_{k,T}) v}{ \nabla w}_{U_k(T)}, \\
  b_{T}(v,w) & = \skal{A(\rchi_T \nabla - \nabla {\QQ }_{k,T}) v}{ \nabla w}_{U_k(T)}, \\
  \tilde b(v,w) & = \sum_{ T \in \TT^{\text{ref}}_H}^{} b_{T}^{\text{ref}}(v,w) + \sum_{T \in {\TT }_H \setminus \TT^{\text{ref}}_H}^{} b_{T}(v,w).
\end{align*}
\item Similarly, define the functional $\tilde c$ for the right hand side correctors by
  \begin{align*}
    c_{T}^{\text{ref}}(v) & = \skal{A_\text{ref} \nabla \RR^{\text{ref}}_{k,T} f_{\text{ref}}}{  \nabla v}_{U_k(T)}, \\
    c_{T}(v) & = \skal{A  \nabla \RR_{k,T}f}{\nabla v}_{U_k(T)}, \\
    \tilde c(v) & = \sum_{ T \in \TT^{\text{ref}}_H}^{} c_{T}^{\text{ref}}(v) + \sum_{T \in {\TT }_H \setminus \TT^{\text{ref}}_H}^{} c_{T}(v),
  \end{align*}
  and solve for $\tilde u^{H}_{k} \in  V_H$ in
  \begin{equation}\label{proposedmethod}
    \begin{aligned}
      \tilde b(\tilde u^{H}_{k},v) = {}& F(v) - \tilde c(v) - \tilde b(g,v) 
    \end{aligned}
  \end{equation}
  for all $v \in V_H$, and compute the solution as
  \begin{equation}\label{eq:approximate_tildeuk}
    \tilde u_{k}= \tilde{u}^H_k -\tilde{Q}_k\tilde{u}^H_k + \tilde{\RR}_k \tilde{f} -\tilde{Q}_k g.
  \end{equation}
\end{enumerate}
\end{definition}

\begin{bemerkung}[Individual marking]
  The described algorithm can be enhanced in terms of efficiency by separating $\TT^{\text{ref}}_H$ for each corrector type $\QQ_{k,T} \lambda_x$, $\RR_{k,T} f$, and $\QQ_{k,T} g$, only updating the correctors with respect to their corresponding error indicator. For example, if $A = A_{\text{ref}}$ but $f \ne f_{\text{ref}}$, obviously the $f$-independent correctors $\QQ_{k,T} \lambda_x$ and $\QQ_{k,T} g$ need not to be recomputed in Step 3, since only $\RR_{k,T} f$ can differ from its reference counterpart $\RR^{\text{ref}}_{k,T} f_{\text{ref}}$. As stated, however, the algorithm would (unnecessarily) recompute the $f$-independent correctors as well. For readability of this paper, we decided to omit individual marking.
\end{bemerkung}

%We remark that the chosen fraction $p$ in point 3\ does most likely correspond to an error tolerance TOL that we want to match such that $e_{u,p} \leq TOL$. Certainly, the method is only efficient for rather small $p$. However, if the perturbation is too big, we want to increase the computational effort instead of getting an inaccurate solution. We will present the role of TOL in Chapter 5 since for now, it is only important that we indeed choose a fraction $p$ regardless of the concrete choice.

%Or formulated in different way.
%\begin{definition}[PG-LOD with adaptively updated correctors]
%	\label{uvc}
%Fix the patch size $k \in \mathbb{N}$ and a certain recompute fraction $p \in [0,1]$. Let TOL(p) be the corresponding tolerance such that
%\begin{equation*}
%	e_u \leq \text{TOL(p)} \qquad \text{ and } \qquad
%		\frac{\# \hat{\QQ }}{\left(\# \QQ  + \# \hat{\QQ }\right)} \approx p,	
%\end{equation*}	
%where $\# \hat{\QQ }$ denotes the number of updated correctors and $\# {\QQ }$ the number of old reference correctors. Let $\tilde V_{k,p}$ be the resulting LOD space with mixed correctors $\hat{\QQ }^T$ and ${\QQ }^T$. The approximation of $u$ in \eqref{eq:perturbed_problem} is given by finding $\tilde u_{k,p} \in \tilde V_{H,k}$ such that, for all $v \in V_H$, it holds that
%\begin{equation}
%	 \tilde a(\tilde u_{k,p},v) = F(v).
%\end{equation}
%\end{definition}

\section{Error analysis}
\label{cha:error_analysis_for}
This section is devoted to the theoretical justification of the proposed method. We present the main theorem of this work. The theorem justifies local recomputation of the correctors based on the value of the error indicators.

%\subsection{Main results} %\tim{keep the sectioning like this?}
%We first of all prove a lemma that justifies the intuition of Definition %\ref{def:indicators}. These indicators were first developed in \cite{2017}.
%\label{sub:main_results}

%In the error analysis we can exploit Lemma \ref{lemma:indicator} to state our the key result for our method.
\begin{theorem}[Error bound for the PG-LOD with adaptively updated correctors] \label{theorem}
  If
  \begin{equation*}
    \max_{T\in \TT^{\text{ref}}_H}(E_{\QQ V_H,T}\left(\norm{f}_{L^2(\Omega)} + \vertiii{g}\right), E_{f,T} + E_{\RR f, T}, E_{\QQ g,T})\leq \text{TOL}\left(\norm{f}_{L^2(\Omega)} + \vertiii{g}\right)
  \end{equation*}
  then there exist $k_0 > 0 $ and $\tau_0 > 0$ such that for all $k > k_0$ and $\tau < \tau_0$ with $\text{TOL} = \tau k^{-d/2}$ the error bound
  \begin{equation} \label{eq:theorem_4.2}
    \vertiii{u - \tilde{u}_k} \lesssim k^{d/2} ( \theta^k  + \text{TOL}) (\norm{f}_{L^2(\Omega)} + \vertiii{g})
  \end{equation}
  is satisfied. Here, $0 < \theta < 1$ is independent of $H$, $k$, $\tau$ and TOL.
\end{theorem}
Before proving the theorem, we require the following lemma.
\begin{lemma}[Error indicators bound the errors in reference correctors and integrals] \label{lemma:indicator}
  For all $v \in V_H$, the following bounds hold
  \begin{align*}
    \vertiii{\QQ _{k,T}v- \QQ ^{\text{ref}}_{k,T}v} &\leq \norm{ A^{-1/2}  ( A-A_{\text{ref}}) (\rchi_T  \nabla v -  \nabla {\QQ ^{\text{ref}}_{k,T}}v)  }_{L^2(U_k(T))} && \hspace{-0.5cm} \leq E_{\QQ V_H,T} \, \vertiii{v}_T,  \\
    \vertiii{\RR _{k,T}f- \RR ^{\text{ref}}_{k,T}f_{\text{ref}}} &\leq E_{f,T} + \norm{ A^{-1/2}  ( A-A_{\text{ref}}) \nabla {\RR ^{\text{ref}}_{k,T}} f_{\text{ref}}  }_{L^2(U_k(T))} && \hspace{-0.5cm} \leq E_{f,T} + E_{\RR f, T}, \\
    \vertiii{\QQ _{k,T}g- \QQ ^{\text{ref}}_{k,T}g} &\leq \norm{ A^{-1/2}  ( A-A_{\text{ref}}) (\rchi_T  \nabla g -  \nabla {\QQ ^{\text{ref}}_{k,T}}g)  }_{L^2(U_k(T))} && \hspace{-0.5cm} \leq  E_{\QQ g,T}.
  \end{align*}
\end{lemma}
\begin{beweis}[Lemma~\ref{lemma:indicator}]
	For any $v \in V_H$, we define $z := \QQ _{k,T}v- \QQ ^{\text{ref}}_{k,T}v \in V^{\text{f}}(U_k(T))$ and observe
	\begin{align*}
		\vertiii{z}_{U_k(T)}^2 &= \skal{ A \nabla v}{ \nabla z}_T
		 - \skal{ A \nabla  \QQ ^{\text{ref}}_{k,T}v  }{ \nabla z}_{U_k(T)} \\
		 & \qquad \qquad  + \skal{ A_{\text{ref}}  \nabla  \QQ ^{\text{ref}}_{k,T}v}{ \nabla z}_{U_k(T)} - \skal{ A_{\text{ref}}  \nabla v}{ \nabla z}_T \\
		 &\leq \norm{ A^{-1/2}  ( A-A_{\text{ref}}) (\rchi_T  \nabla v -  \nabla {\QQ ^{\text{ref}}_{k,T}}v)  }_{L^2(U_k(T))} \cdot \vertiii{z}_{U_k(T)},
	\end{align*}
        which yields the first inequality of the first part. We proceed to get the second inequality,
\begin{align*}
		\vertiii{z}_{U_k(T)}^2 &\leq \norm{ A^{-1/2}  ( A-A_{\text{ref}}) (\rchi_T  \nabla v -  \nabla {\QQ ^{\text{ref}}_{k,T}}v)  }^2_{L^2(U_k(T))} \\
&\leq \max_{ \substack{w|_{T},\,w \in V_H %, \\ \|A_\text{ref}^{1/2}\nabla w\|_{L^2(T)}=1
}} \frac{\norm{ A^{-1/2}  ( A-A_{\text{ref}}) (\rchi_T  \nabla w -  \nabla {\QQ ^{\text{ref}}_{k,T}}w)  }^2_{L^2(U_k(T))}}{\|A_\text{ref}^{1/2}\nabla w\|^2_{L^2(T)}}\|A_\text{ref}^{1/2}\nabla v\|^2_{L^2(T)}\\
&\leq \|A_\text{ref}A^{-1}\|_{L^\infty(T)} \mkern-15mu \sum_{\substack{T' \in \TT _H \\ T' \cap U_k(T) \neq 0}}^{}\mkern-15mu \norm{\delta}^2_{L^{\infty}(T')}\mkern-6mu \max_{\substack{ w |_{T}, \,w \in V_H%,\\ \|A_\text{ref}^{1/2}\nabla w\|_{L^2(T)}=1
}}\mkern-12mu \frac{\norm{ A_{\text{ref}}^{1/2} (\rchi_T  \nabla w -  \nabla {\QQ ^{\text{ref}}_{k,T}}w)  }^2_{L^2(T')}}{\|A_\text{ref}^{1/2}\nabla w\|^2_{L^2(T)}}\vertiii{v}^2_T\\
&= E_{\QQ V_H,T}^2\vertiii{v}^2_T.
	\end{align*}
	with $\delta = A^{-1/2}(A-A_{\text{ref}}) A_{\text{ref}}^{-1/2}$.

	The second result follows analogously with $z':= \RR _{k,T}f - \RR ^{\text{ref}}_{k,T}f_{\text{ref}} \in V^{\text{f}}(U_k(T))$ and
	\begin{align*}
		\vertiii{z'}^2_{U_k(T)} & = \skal{f-f_{\text{ref}}}{z'}_T + \skal{f_{\text{ref}}}{z'}_T - \skal{A  \nabla \RR^{\text{ref}}_{k,T}f_{\text{ref}}}{ \nabla z'}_{U_k(T)} \\
                                        &\le \Big( \|A^{-1}\|^{1/2}_{L^\infty(U(T))} \norm{f - f_{\text{ref}}}_{L^2(T)} \sup_{\substack{v \in V}} \frac{\|(\id-\mathcal{I}_H)v\|_{L^2(T)}}{\|\nabla v\|_{L^2(U(T))}} + {}\\
                                        &\phantom{{}=\Big(} \norm{A^{-1/2}(A - A_{\text{ref}}) \nabla  \RR^{\text{ref}}_{k,T}f_{\text{ref}}}_{L^2(U_k(T))} \Big) \cdot \vertiii{z'}_{U_k(T)}
	\end{align*}        
	Similar arguments yield the third result of the lemma.
\end{beweis}

\begin{beweis}[Theorem~\ref{theorem}]
  The full error is $\vertiii{u - \tilde u_k} \le \vertiii{u - u_k} + \vertiii{u_k - \tilde u_k}$. The first term stems from the localization and is bounded by $\vertiii{u - u_k} \lesssim k^{d/2}\, \theta^k \left(\|f\|_{L^2(\Omega)} + \vertiii{g} \right)$ according to \eqref{error_bound_1}.

  Before proceeding with the second term, we note that we can bound the error in the global correctors $\tilde{\QQ}_k$ and $\tilde{\RR}_k$ in terms of the patch overlap $k^{d/2}$ and TOL by using Lemma~\ref{lemma:indicator} and the assumption on \text{TOL} stated in this theorem: for all $v \in V_H$,
  \begin{equation} \label{splitting}
    \begin{split}
      \vertiii{({\QQ }_k - \tilde{\QQ}_k)v}^2 & = \vertiii{\sum\nolimits_{T \in \TT^{\text{ref}}_H}^{} ({\QQ }_{k,T}- {\QQ }^\text{ref}_{k,T})v}^2 \\
      & \lesssim \sum\nolimits_{T \in \TT^{\text{ref}}_H}^{}k^d \vertiii{({\QQ }_{k,T}- {\QQ }^\text{ref}_{k,T})v}^2\\			
      &\lesssim k^d \text{TOL}^2 \vertiii{v}^2.
    \end{split}
  \end{equation} 
  Analogously, we get $\vertiii{(\RR_k - \tilde{\RR}_k)f} + \vertiii{(\QQ_k - \tilde{\QQ}_k)g} \lesssim k^{d/2}\text{TOL}(\|f\|_{L^2(\Omega)} + \vertiii{g})$. Using \eqref{eq:decay} and \eqref{splitting}, we additionally bound
  \begin{equation*}
    \vertiii{\tilde{\QQ}_k v} = \vertiii{\QQ v} + \vertiii{(\QQ - \QQ_k) v} + \vertiii{(\QQ_k - \tilde{\QQ}_k) v} \lesssim (1 + k^{d/2} \theta^k + k^{d/2} \text{TOL}) \vertiii{v} \lesssim \vertiii{v}.
  \end{equation*}
  Next, we proceed with the second term, using \eqref{eq:localized_uk} and \eqref{eq:approximate_tildeuk} and the bounds above,
  \begin{equation}\label{eq:error_in_tilde}
    \begin{aligned}
      \vertiii{u_k - \tilde u_k} & \le \vertiii{u^{H}_k - \tilde u^{H}_k} + \vertiii{\QQ_k u^H_k - \tilde \QQ_k \tilde u^H_k} + \vertiii{\RR_kf - \tilde \RR_k \tilde f} - \vertiii{\QQ_k g - \tilde \QQ_k g} \\
      & \lesssim \vertiii{u^{H}_k - \tilde u^{H}_k} + \vertiii{\QQ_k u^H_k - \tilde \QQ_k u^H_k} + \vertiii{\RR_kf - \tilde \RR_k \tilde f} - \vertiii{\QQ_k g - \tilde \QQ_k g} \\
      & \lesssim \vertiii{u^{H}_k - \tilde u^{H}_k} + k^{d/2} \text{TOL}(\|f\|_{L^2(\Omega)} + \vertiii{g}),
    \end{aligned}
  \end{equation}
  where we added $\pm \tilde \QQ_k u^H_k$ in the second estimate and use $\vertiii{u^H_k} \mkern-3mu \lesssim \mkern-3mu \|f\|_{L^2(\Omega)} + \vertiii{g}$ for the last step.

  It remains to bound the energy norm of $u^{H}_k - \tilde u^{H}_k \in V_H$. The first step is to establish a coercivity inequality for $\tilde b$ on $V_H$. To do this, we define the auxiliary bilinear form $b = \sum_{T \in \TT_H} b_T$ (with $b_T$ defined in Definition~\ref{procedure}) and bound the consistency error $[b - \tilde b](v, w)$ for $v, w \in V_H$,
  \begin{align}\label{eq:consistency_error}
    [b - \tilde b](v, w) = \sum_{T \in \TT^{\text{ref}}_H} [b^{\text{ref}}_T - b_T](v, w) \lesssim k^{d/2} \text{TOL} \vertiii{v} \vertiii{w},
  \end{align} 
  \vspace{-0.1cm}
  since
  \begin{align*}
    [b^{\text{ref}}_T - b_T](v, w) & = \skal{A_\text{ref}(\rchi_T  \nabla  -  \nabla \QQ ^{\text{ref}}_{k,T}) v - A(\rchi_T \nabla - \nabla {\QQ }_{k,T}) v}{ \nabla w}_{U_k(T)} \\
                                       & = \skal{(A_\text{ref} - A)(\rchi_T \nabla  -  \nabla \QQ ^{\text{ref}}_{k,T}) v - A\nabla ({\QQ }_{k,T} - \QQ ^{\text{ref}}_{k,T}) v}{ \nabla w}_{U_k(T)} \\
                                       & \le \norm{A^{-1/2}(A_\text{ref} - A)(\rchi_T \nabla  -  \nabla \QQ ^{\text{ref}}_{k,T})v}_{L^2(U_k(T))}\vertiii{w}_{U_k(T)}\\
                                       & \phantom{\le{}} + \vertiii{({\QQ}_{k,T} - \QQ ^{\text{ref}}_{k,T}) v} \vertiii{w}_{U_k(T)} \\
                                       & \lesssim \text{TOL} \vertiii{v}_T \vertiii{w}_{U_k(T)},
  \end{align*}
  where we again use Lemma~\ref{lemma:indicator}. Additionally, for $v \in V_H$, we note that
  \begin{equation}\label{eq:pg_coercivity}
    a(v - \QQ v, v) = a(v - \QQ v, v - \QQ v) = \vertiii{v - \QQ v}^2 \ge \alpha\beta^{-1}C_{\mathcal{I}_H}^{-2} \vertiii{\mathcal{I}_H (v - \QQ v)}^2 \gtrsim \vertiii{v}^2,
  \end{equation}
  and further that $a(v - \QQ_kv, v) = b(v,v)$. Using \eqref{eq:pg_coercivity}, \eqref{eq:decay}, \eqref{eq:consistency_error}, and the triangle inequality, we get for $v \in V_H$,
  \begin{equation*}
    \begin{aligned}
      |\tilde b(v, v)| & \ge a(v - \QQ v, v) - |a((\QQ - \QQ_k)v, v)| - |[b - \tilde b](v, v)| \\
      & \ge (C_1 - C_2 k^{d/2} \theta^{k} - C_3 k^{d/2} \text{TOL} )\vertiii{v}^2,
    \end{aligned}
  \end{equation*}
  where $C_1$, $C_2$, and $C_3$ are independent of $H$, $k$ and TOL. From this inequality, we note that there exist $k_0$ and $\tau_0$ such that, for all $k > k_0$ and $\tau < \tau_0$ (with $\text{TOL} = \tau k^{-d/2}$), we have the coercivity inequality
  \begin{equation*}
    \begin{aligned}
      \tilde b(v, v) \ge \tilde \gamma \vertiii{v}^2 \gtrsim \vertiii{v}^2,
    \end{aligned}
  \end{equation*}
  for all $v \in V_H$, where $\tilde \gamma$ depends on $k_0$ and $\tau_0$ but not on $H$, $k$, $\tau$ or TOL.
  
  We define $c = \sum_{T \in \TT_H} c_T$ and observe that the localized problem \eqref{eq:PGLODref} can be
  expressed in the same manner as the proposed approximation in \eqref{proposedmethod} as
  \begin{equation*}
    b(u^H_k, v) = F(v) - c(v) - b(g, v),
  \end{equation*}
  for all $v \in V_H$.
  Next, we use the coercivity inequality with $v = u^{H}_k - \tilde u^{H}_k =: e \in V_H$ and the reformulation of \eqref{eq:PGLODref} together with \eqref{proposedmethod} to obtain
  \begin{equation}\label{eq:error_in_vh}
    \begin{aligned}
      \vertiii{u^{H}_k - \tilde u^{H}_k}^2 & \lesssim \tilde b(u^{H}_k - \tilde u^{H}_k, e) = \tilde b(u^{H}_k, e) - F(e) + \tilde c(e) + \tilde b(g,e) \\
      & = [\tilde b - b](u^H_k, e) + [\tilde c - c](e) + [\tilde b - b](g, e) \\
      & \lesssim k^{d/2} \text{TOL} (\vertiii{u^H_k} + \|f\|_{L^2(\Omega)} + \vertiii{g})\vertiii{e} \\
      & \lesssim k^{d/2} \text{TOL} (\|f\|_{L^2(\Omega)} + \vertiii{g})\vertiii{u^{H}_k - \tilde u^{H}_k}, \\
    \end{aligned}
  \end{equation}
  where we used \eqref{eq:consistency_error} for $[\tilde b - b](u^H_k, e)$ (and analogously derived results for $[\tilde c - c](e)$ and $[\tilde b - b](g, e)$).
  Combining the last bound with \eqref{error_bound_1} and \eqref{eq:error_in_tilde} yields the asserted bound.
\end{beweis}

\begin{bemerkung}[Right hand side correction]
	It is possible to ignore $\RR$ in the method. This leads to an additional error term that is proportional to $H$ in \eqref{eq:theorem_4.2}. However, for localized right hand sides or if high accuracy is needed $\RR$ should be included, see \cite{HM14} where this term was first analyzed for LOD.   
\end{bemerkung}

\begin{bemerkung}[$V=H^1_0(\Omega)$]\label{fullydisc}
  The assumption that $V = V_h$ (a finite element space) means that the error will be with respect to the finite element approximation $u\in V_h$. The same analysis goes through if we instead let $V=H^1_0(\Omega)$ but the corresponding continuous solution is not computable. It is assumed that there is a fine mesh with mesh size $h$ for which the error in the fine scale finite element solution is small enough. An additional a priori error bound then gives the full error with respect to the continuous solution using the triangle inequality.
\end{bemerkung}

\section{Implementation}
\label{sec:implementation}
This section discusses implementation details specific to the
presented method, with emphasis on the computation of the error
indicators, parallel computations and memory consumption for
large-scale problems. For implementation details on the LOD
corrector problems we refer to \cite{EPMH16}. However, we would like
to emphasize that the localized computations of correctors on patches
makes it possible to avoid any global computations in the (typically
very large) space $V_h$.

\subsection{Computing the error indicators}
It is important for the method that the local error indicators $E_{\QQ V_H,T}$, $E_{f,T}$, $E_{\RR f,T}$, and $E_{\QQ g,T}$ can be computed efficiently.
Consider the definition of $E_{\QQ V_H,T}$ from \eqref{indidef}:
\begin{equation*}
  E^2_{\QQ V_H,T} := \|A_\text{ref}A^{-1}\|_{L^\infty(T)}\sum_{\substack{T' \in \TT _H \\ T' \cap U_k(T) \neq 0}}^{} \mkern-10mu \norm{\delta}^2_{L^{\infty}(T')}
  \underbrace{\max_{\substack{ w |_{T}, \,w \in V_H}} \frac{\norm{ A_{\text{ref}}^{1/2} (\rchi_T  \nabla w -  \nabla {\QQ ^{\text{ref}}_{k,T}}w)  }^2_{L^2(T')}}{\|A_\text{ref}^{1/2}\nabla w\|^2_{L^2(T)}}}_{{}=: \tilde \mu_{T,T'}}.
\end{equation*}
We denote the maximum factor by $\tilde \mu_{T,T'}$ and note that computing the error corrector for an element $T$ amounts to computing $\|A_\text{ref}A^{-1}\|_{L^\infty(T)}$, and additionally $\norm{\delta}^2_{L^{\infty}(T')} = \norm{A^{-1/2}(A_{\text{ref}}-A) A_{\text{ref}}^{-1/2}}^2_{L^{\infty}(T')}$ and $\tilde \mu_{T,T'}$ for the coarse elements $T'$ in the patch $U_k(T)$. The values $\tilde \mu_{T,T'}$ can be precomputed as $\tilde \mu_{T, T'} = \max_l \mu_l$, where $\mu_l$ is the solution of the eigenvalue problem
\begin{equation} \label{ls_correctors}
	\mathcal{B} x_l = \mu_l \mathcal{C} x_l,
\end{equation}with
\begin{equation} \label{BandC}
	\begin{split}
			\mathcal{B}_{ij} &= \skal*{(  A_{\text{ref}} ( \rchi_T  \nabla \lambda_j -  \nabla \mathcal{Q}_{k,T}^{\text{ref}} \lambda_j )}{ \rchi_T  \nabla \lambda_i -  \nabla \mathcal{Q}_{k,T}^{\text{ref}} \lambda_i }_{T'}, \\
			\mathcal{C}_{ij} &= \skal*{A_{\text{ref}}  \nabla \lambda_j}{ \nabla \lambda_i}_T,
	\end{split}
\end{equation}
for $i,j = 1, \dots, m-1$ and where $m$ denotes the number of basis functions with support in element $T$. For 2D quadrilateral mesh elements this results in a $3 \times 3$ system. 
Given a perturbed coefficient $A$, we then have to compute only the $\norm{\delta}^2_{L^{\infty}(T')}$ terms and multiply them with the precomputed $\tilde \mu_{T,T'}$ before summing, and also computing the $\norm{A_{\text{ref}}A^{-1}}_{L^{\infty}(T)}$ factor. An important consequence of this procedure is that the fine scale reference correctors $\mathcal{Q}_{k,T}^{\text{ref}} \lambda_j$ can be discarded and only the coarse scale quantity $\tilde \mu_{T,T'}$ needs to be stored.

The error indicator $E_{\QQ g,T}$ for ${\QQ}^{\text{ref}}_{k,T}g$ can be computed in a similar way as $E_{\QQ V_H,T}$, the difference being that no eigenvalue problem over functions in $V_H$ restricted to $T$ needs to be solved since $g$ is known a-priori.

The right hand side error indicator $E_{f,T} + E_{\RR f,T}$ is slightly different, since perturbations of both $f$ and $A$ affect the accuracy of ${\RR }^{\text{ref}}_{k,T} f_{\text{ref}}$. The second term $E_{\RR f,T}$ can be computed similarly to $E_{\QQ V_H,T}$ and $E_{\QQ g,T}$. The first term $E_{f,T}$ contains an element local Poincaré-type inequality for functions in the fine space (i.e.\ the null space of $\mathcal{I}_H$),
\begin{equation}
  \label{eq:supremum}
  \nu_T := \sup_{\substack{v \in V}} \frac{\|(\id-\mathcal{I}_H)v\|_{L^2(T)}}{\|\nabla v\|_{L^2(U(T))}}.
\end{equation}
We have from the a-priori approximability bound of $\mathcal{I}_H$ in
\eqref{Ia} that this quantity scales with $H$. If the constant
$C_{\mathcal{I}_H}$ is known it can be used here. In many practical
situations it is possible to get a sharp estimate by computational
means since $V$ is the finite element space $V_h$ and the
bound we seek is the maximum eigenvalue to a problem posed on a
1-layer element patch $U(T)$ on the fine scale. For all the
experiments in Section~\ref{sec:experiments} we discretize the two
dimensional unit square with a uniform grid and get three possible
cases: $U(T)$ consists 4, 6 or 9 elements (depending on whether it is
in the corner, on the edge or in the interior of the domain). By
computing \eqref{eq:supremum} with $V = V_h$ for a
decreasing range of small $h$, it was possible to obtain the estimate
$\nu_T \le 0.25H$ for all $T$. This value was used in the experiments
in Section~\ref{sec:experiments}.

\subsection{Algorithm, parallelization and memory consumption}
\begin{algorithm}
	\SetAlgoLined
        \SetKwInOut{Input}{Input}
        \SetKwInOut{Output}{Output}
        \Input{Reference data ($A_{\text{ref}}$, $f_{\text{ref}}$)}
	Pick $k$ and TOL \\
        \For{\textup{all} $T$}{
          Precompute ${\QQ }^{\text{ref}}_{k,T} \lambda_j$, ${\RR }^{\text{ref}}_{k,T} f_{\text{ref}}$ and ${\QQ}^{\text{ref}}_{k,T}g$ for all $j$ (discard at end of iteration)\\
          Precompute and save $\nu_T$ and $\tilde \mu_{T,T'}$ for all $T' \subset U_k(T)$ \\
          Precompute and save $b^{\text{ref}}_T(\lambda_j, \lambda_i)$, $c^{\text{ref}}_T(\lambda_i)$, and $b^{\text{ref}}_T(g, \lambda_i)$ for all $i$ and $j$ \\
        }
        \For{\textup{perturbed data} ($A$, $f$)}{
          ${\TT}^{\text{ref}}_H \leftarrow \emptyset$ \\
          \For{\textup{all} $T$}{
            Compute $E_{\QQ V_H,T}$, $E_{f, T}$, $E_{\RR f,T}$ and $E_{\QQ g,T}$ using $\tilde \mu_{T,T'}$ and $\nu_T$ \\
            \uIf{$\max(E_{\QQ V_H,T}(\norm{f}_{L^2(T)} + \vertiii{g}), E_{f, T}, E_{\RR f,T}, E_{\QQ g,T}) \leq \textup{TOL}$}{
              ${\TT}^{\text{ref}}_H \leftarrow \{T\} \cup {\TT}^{\text{ref}}_H$ \\
              Compute ${\QQ }_{k,T} \lambda_j$ for all $j$, ${\RR }_{k,T}f$, and ${\QQ }_{k,T}g$ (discard at end of iteration) \\
              Compute and save $b_T(\lambda_j, \lambda_i)$, $c_T(\lambda_i)$, and $b_T(g, \lambda_i)$ for all $i$ and $j$ \\
            }
          }
          Assemble stiffness matrix $\tilde{K}_{i,j} = \sum_{T \in {\TT}^{\text{ref}}_H} b^{\text{ref}}_T(\lambda_j, \lambda_i) + \sum_{T \notin {\TT}^{\text{ref}}_H} b_T(\lambda_j, \lambda_i)$ and load vector similarly \\
          Solve for $\tilde u_k^H$ using \eqref{proposedmethod} \\
          Possibly compute full solution by $\tilde u_k = (\id - \tilde{\QQ}_k)\tilde u_k^H - \tilde{\QQ}_k g + \tilde{\RR}_k \tilde f$
        }
        \caption{Simplified procedure of the proposed method}
        \label{pseudo_code}
\end{algorithm}
Algorithm~\ref{pseudo_code} shows an example of how the computational steps can be carried out.
It is based on the assumption that we have to solve for multiple perturbations of a certain reference coefficient and therefore consists of an initial stage when reference quantities are computed. We make the observations that:
\begin{itemize}
\item The amount of data stored from the initial stage scales like $k^{d}H^{-d}$ since it consists only of overlapping quantities on the coarse mesh.
\item The initial stage can be computed in parallel over $T$ thanks to the PG-LOD formulation since  $b^{\text{ref}}_T(\lambda_j, \lambda_i)$, $c^{\text{ref}}_T(\lambda_i)$, and $b^{\text{ref}}_T(g, \lambda_i)$ depend only on the correctors for $T$.
\end{itemize}

Following the initial stage, a loop over the perturbed coefficients
follows. For each perturbed coefficient, there is another loop over
the elements $T$. It computes updated correctors for the elements for
which the error indicators tell it is needed. Finally, the coarse
scale linear system is assembled and solved and the full solution
can be computed. We note that:
\begin{itemize}
  \item The iterations in the loop over the perturbed coefficients can be executed in parallel.
  \item The iterations in the loop over the elements $T$ can be executed in parallel.
  \item In order to compute the error indicators, the reference coefficient $A_{\text{ref}}$ needs to either be stored explictly (amount $h^{-d}$), be transmitted patch-wise (amount $k^{d}(h/H)^{-d}$) between computers, or be generated from a low-dimensional representation on demand.
  \item There is a reduction over $T$ to assemble the stiffness matrix and the load vector, but only coarse scale data of amount $k^{d}H^{-d}$ is needed for this reduction. This means that fine scale information does not need to be transmitted between computers during reduction.
  \item The coarse scale solution $\tilde u^H_k$ is readily available after solving the coarse scale linear system. If the full solution $\tilde u_k$ is requested in some area of the domain, the correctors in that area have to be recomputed or stored from a previous corrector computation.
  \end{itemize}

\begin{bemerkung}[Periodicity]
  Applications such as composite materials often lead to a periodic structure of the underlying reference coefficient. 
  In this case, correctors can be reused. 
  In the full periodic case this leads to only one corrector problem for full patches and comparably few for the boundary patches. 
  This means that memory consumption as well as complexity decreases significantly. 
\end{bemerkung}
  
\section{Numerical experiments}
\label{sec:experiments}
In this section we present three experiments where the diffusion coefficient is perturbed by defects with either no domain mapping, 
local domain mapping or global domain mapping. We performed this experiments with \texttt{gridlod} in Python (see \cite{gridlod}). The entire code for all our experiments is available on GitHub in \cite{gridlod-on-perturbations}.
Our experiments are performed on a $2$D quadrilateral mesh on $\Omega = [0,1]^2$. For the fine and the coarse scale discretization $V_h$ and $V_H$, 
we use standard $\mathcal{P}_1$ finite element spaces on a fine mesh $\mathcal{T}_h$ with $256 \times 256$ fine elements and a coarse mesh $\mathcal{T}_H$ with $32 \times 32$ coarse elements. 
We use $V_h$ for the fine reference solution and for the fine scale discretization of each corrector problem. 
For the mesh size $H = 1/32$ we conclude from the error estimate \eqref{error_bound_1} that the localization parameter $k = 4 \geq \abs{\log(H)}$ suffices. 

\begin{figure}
	\begin{subfigure}[b]{0.3\textwidth}
		\includegraphics[width=\textwidth]{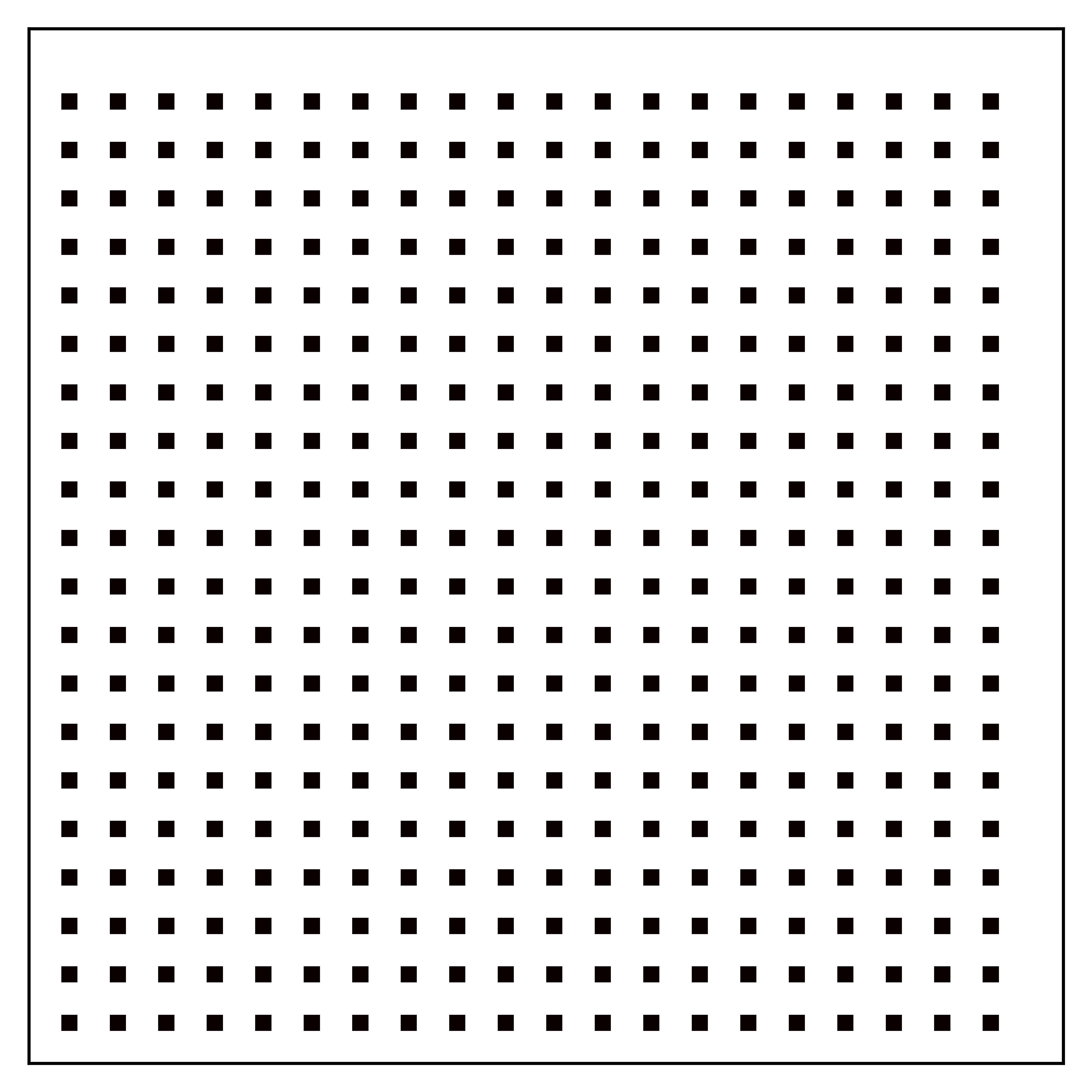}
	\end{subfigure}
	\begin{subfigure}[b]{0.3\textwidth}
		\includegraphics[width=\textwidth]{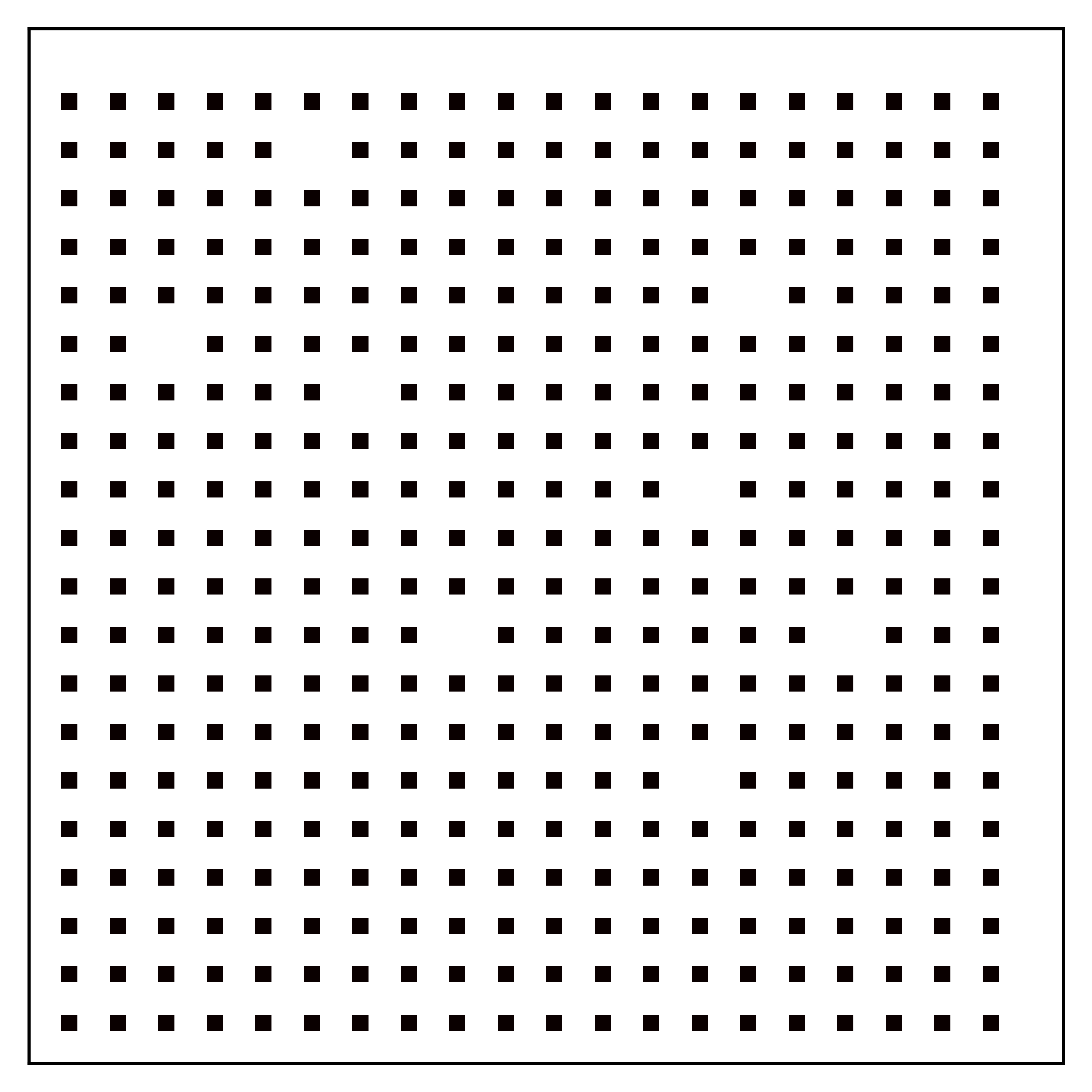}
	\end{subfigure}
	\centering
	\caption{Reference coefficient $A_{\text{ref}}$ (left) and defect perturbation $A$ (right). Black is $1$, white is $0.1$.}
	\label{ex1:coefficients}
\end{figure}

In all our experiments we use a reference coefficient that is piece-wise constant on every fine mesh element. 
Furthermore $A_{\text{ref}}$ can be expressed as stated in Section \ref{sec:perturbations}, i.e.\ it takes two values $1$ and $\alpha$.
For the background $\Omega_{\alpha}$ we choose $\alpha=0.1$. 
In addition, $\Omega_1$ is always a non connected subdomain of $\Omega$ which is conforming with respect to $\mathcal{T}_h$. 
For the sake of convenience we neglect an explicit definition of $A_{\text{ref}}$ and $A$ as we visualize them in the figures. 
We call a perturbation a (local) defect when a (fully connected) subdomain of $\Omega_1$ becomes a subdomain of $\Omega_{\alpha}$. 
In all our experiments $\Omega_1$ is represented by black squares and a defect means that a square gets equalized to the white background (compare Figure \ref{ex1:coefficients}). 
These defects always occur with a probability of $2\%$. The right hand side $f_y$ is always defined by $f_y(y) = \chi_{[1/8,7/8]^2}(y)$. 
For the sake of simplicity we use zero Dirichlet boundary conditions $g \equiv 0$. 

In the experiments we consider the relative error
$$
\mathcal{E}_{\text{rel}}(\tilde{u}_k,u_k) = \frac{\vertiii{\tilde{u}_k - u_k}}{\vertiii{\tilde{u}_k}},
$$ 
where $u_k$ is the best PGLOD solution of \eqref{eq:localized_uk} and $\tilde{u}_k$ is the solution of \eqref{eq:approximate_tildeuk} for a specific the tolerance TOL in the algorithm of Definition \ref{procedure}. 
For $\text{TOL}=\infty$, we clearly have $0\%$ updates of the correctors whereas $\text{TOL}= 0$ corresponds to $100\%$ updates. 
For $100\%$ updates we then end up with the standard Petrov--Galerkin LOD error that is dependent on our data and discretizations. 
In order to observe the complete behavior of $\mathcal{E}_{\text{rel}}(u_k,\tilde{u}_k)$, 
we compute $\mathcal{E}_{\text{rel}}(u_k,\tilde{u}_k)$ for every possible choice of TOL (and thus for every percentage of updates).
The relative best PGLOD error $\mathcal{E}_{\text{rel}}(u,u_k)$ is always around $10^{-3}$ which means that we are comparing to a sufficiently accurate solution. 
\subsection{Defects}
\begin{figure}[h]
	\begin{subfigure}[b]{0.32\textwidth}
		\includegraphics[width=\textwidth]{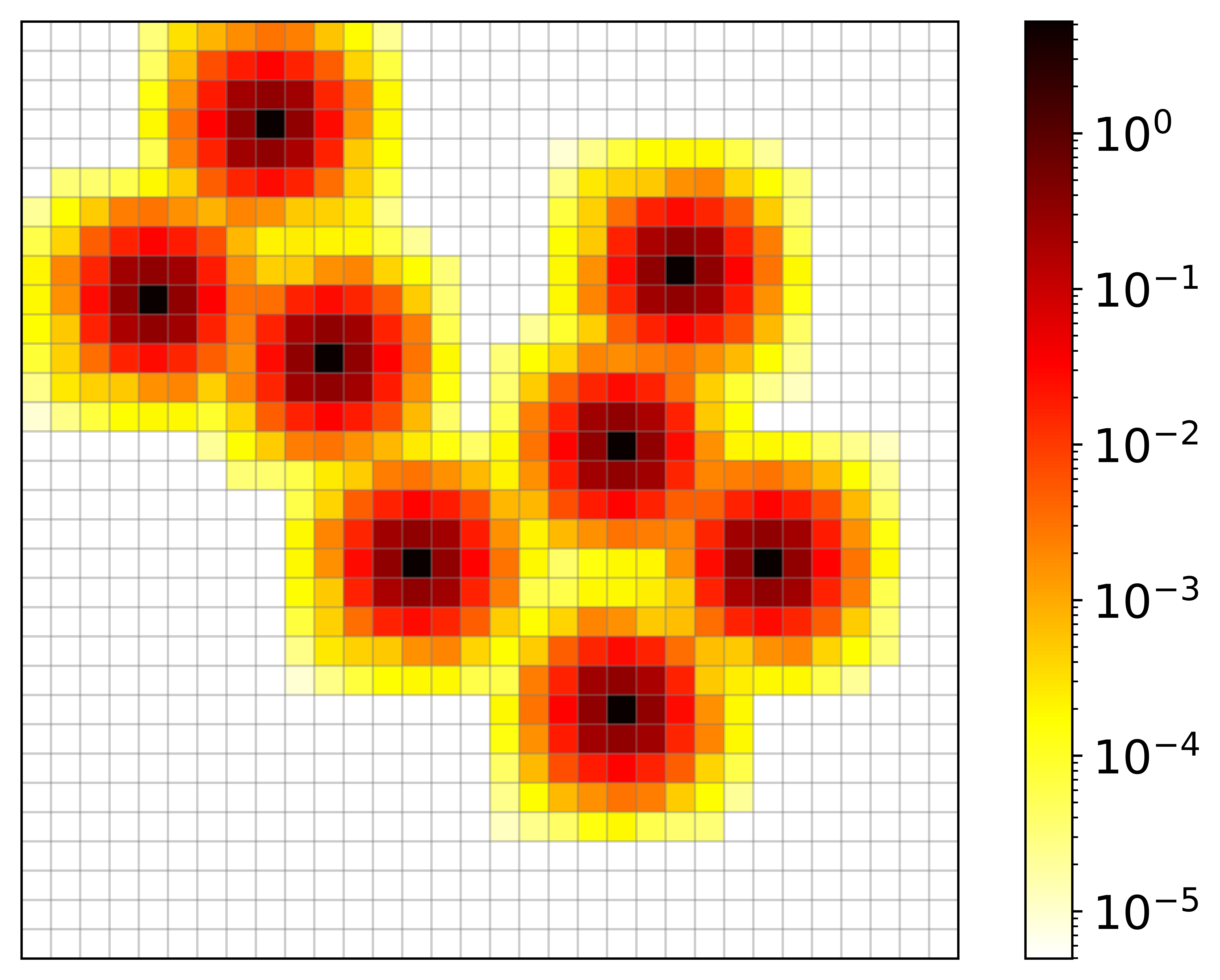}
	\end{subfigure}
	\begin{subfigure}[b]{0.32\textwidth}
		\includegraphics[width=\textwidth]{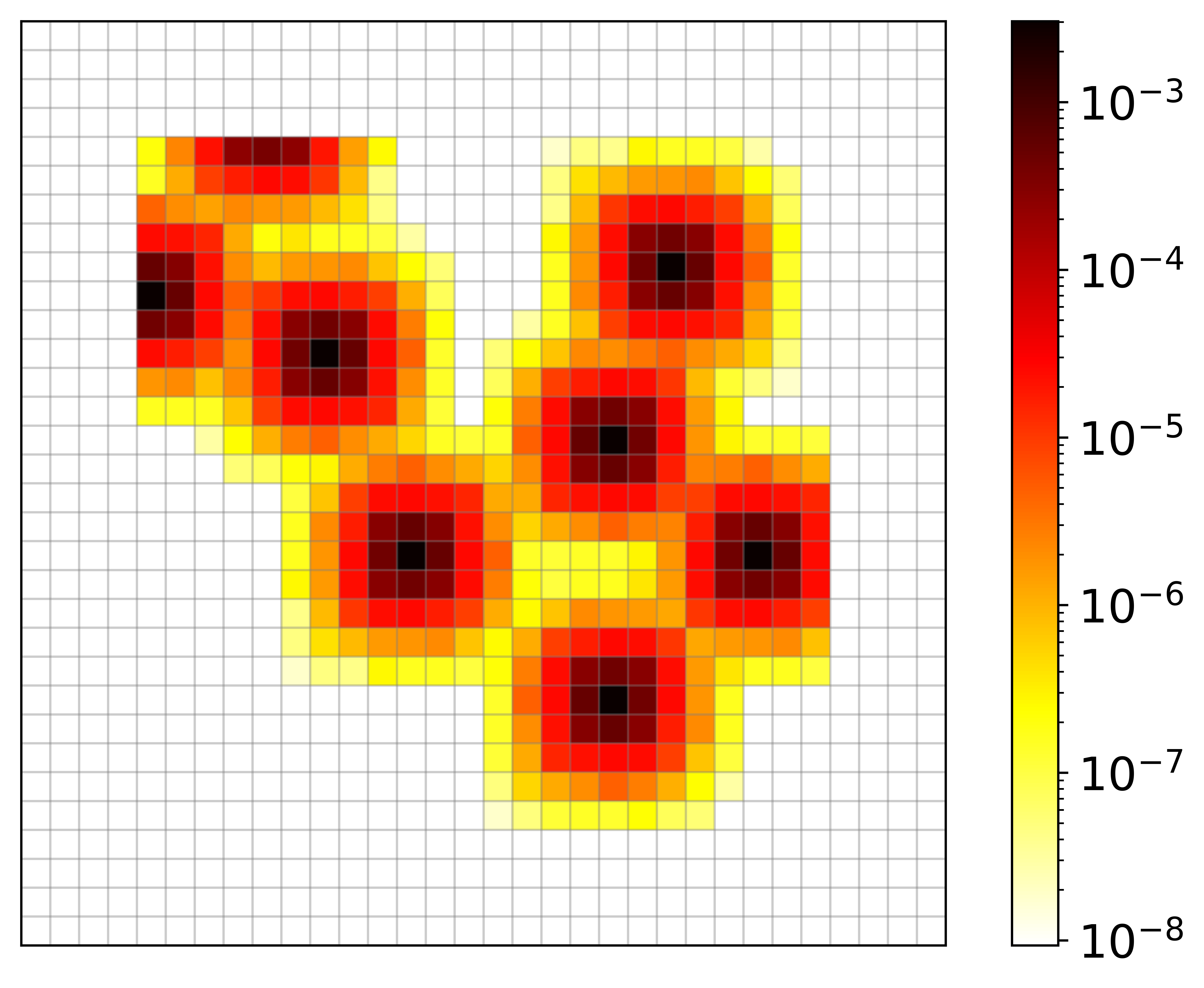}
	\end{subfigure}
	\centering
	\caption{Error indicators $E_{\QQ V_H,T} \cdot \norm{f}_{L^2(\Omega)}$  (left) and ${E_{\RR f,T}}$ (right).}\label{ex1:indicators}
\end{figure}

\begin{figure}[h]
	\centering
	\includegraphics[width=0.55\textwidth]{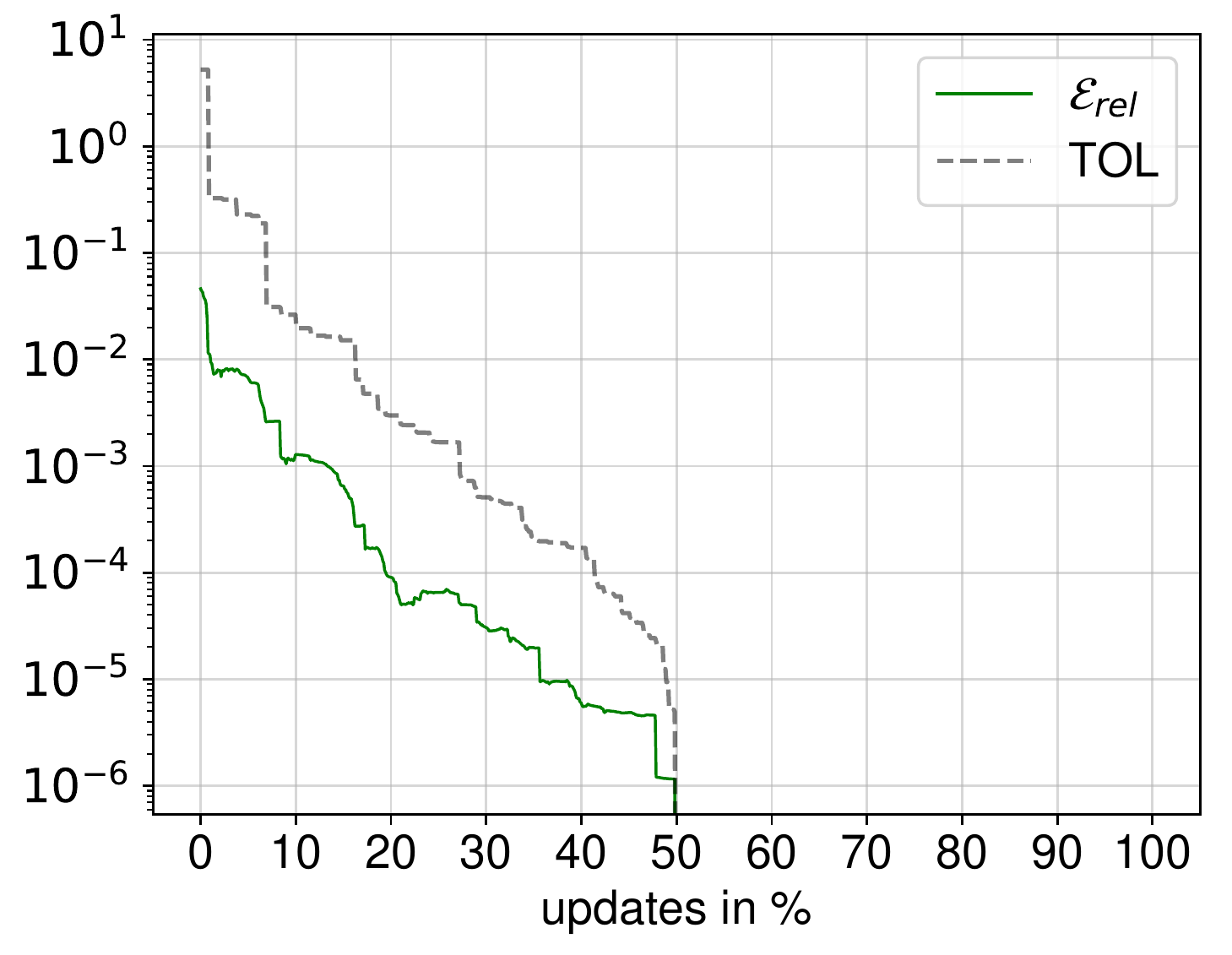}
	\caption{Relative error improvement for defects.}
	\label{ex1:errors}
\end{figure}

In the first experiment we let $\psi=\id$ which means we only consider defects. 
Figure \ref{ex1:coefficients} displays the coefficient and its perturbation whereas the error indicators $E_{\QQ V_H,T}$ and ${E_{\RR f,T}}$ are plotted for each $T$ in Figure \ref{ex1:indicators}. Note that ${E_{f,T}}$ and ${E_{\QQ g,T}}$ are zero for this example. 
The coarse mesh is visible in the background. We can clearly see that the error indicator $E_{\QQ V_H,T}$ detects the defects in the coefficient correctly. 
Furthermore $E_{\QQ V_H,T}$ is exponentially decaying away from each defect. 
In a view of ${E_{\RR f,T}}$, we see that its support coincided with the support of $f_{\text{ref}}$ since for $\psi = \id$ we have $f = f_{\text{ref}}$ which means that whenever $f$ is zero on $T$, 
both $\RR_{k,T}$ and $\RR^{\text{ref}}_{k,T}$ are zero and thus $E_{\RR f,T} = 0$. 
We also observe that $E_{\QQ V_H,T}$ is significantly greater than $E_{\RR f,T}$.
%\axel{My only concern is that we do not say much about what tol is for different levels of updates. Could we pick a resonable value TOL=0.05 (I dont know the value but you choose it so you get roughly 8\% in the first example and a nice number) and then a smaller one which corresponds to say 20\% in the first example so the reader gets an idea of the size of TOL? It is also important since its in the name of the error $e_TOL$. In the best of worlds the same values of TOL also makes sense in the other two examples...}
%\tim{I added TOL to the plot}
From $\mathcal{E}_{\text{rel}}(u_k,\tilde{u}_k)$ in Figure \ref{ex1:errors} we see a high improvement for few updates of the correctors and a sufficiently fast convergence to the best PGLOD solution. 
From this experiment we conclude that the method can efficiently be used for local defects.


\subsection{Local domain mappings}
\begin{figure}[h]
	\begin{subfigure}[b]{0.3\textwidth}
		\includegraphics[width=\textwidth]{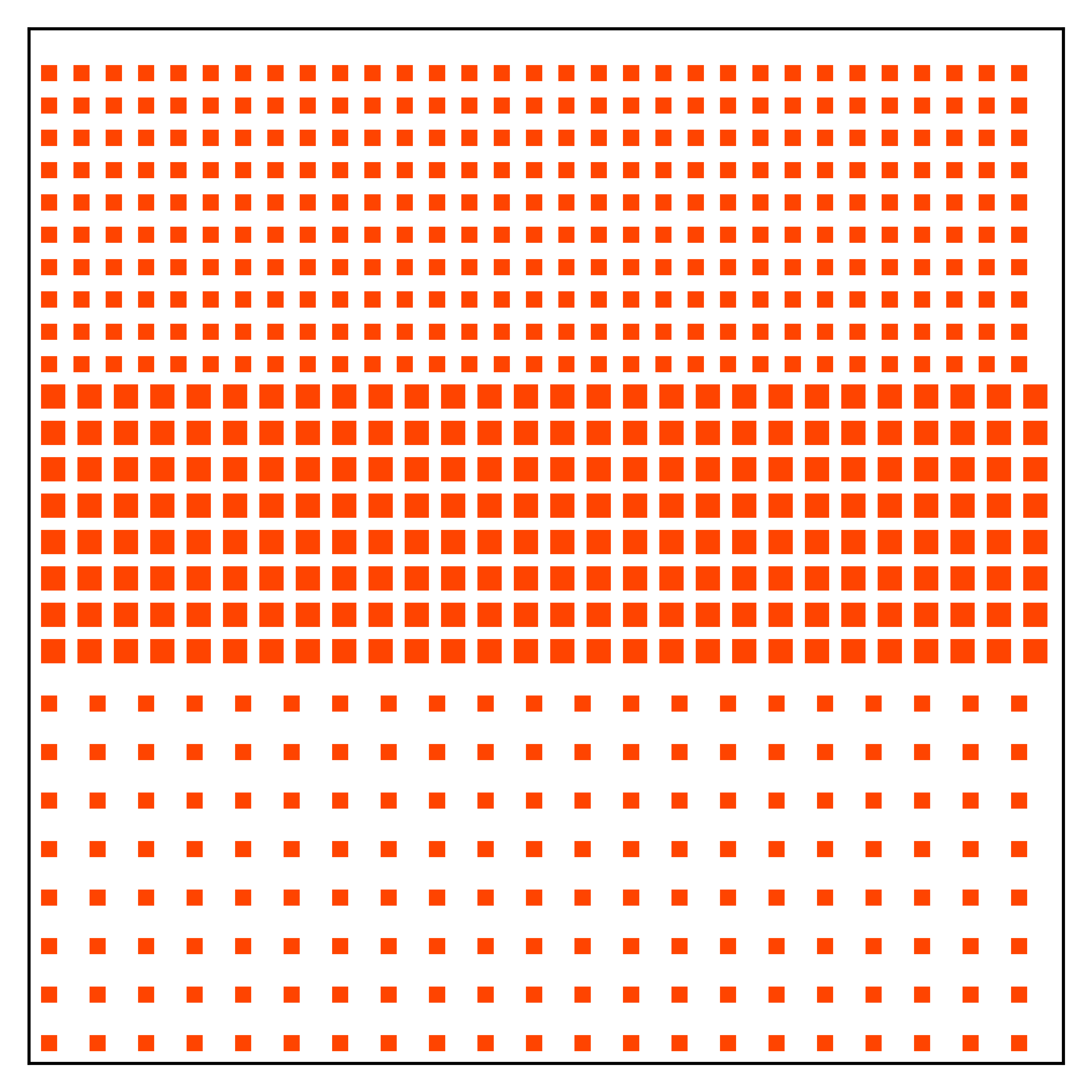}
	\end{subfigure}
	\begin{subfigure}[b]{0.3\textwidth}
		\includegraphics[width=\textwidth]{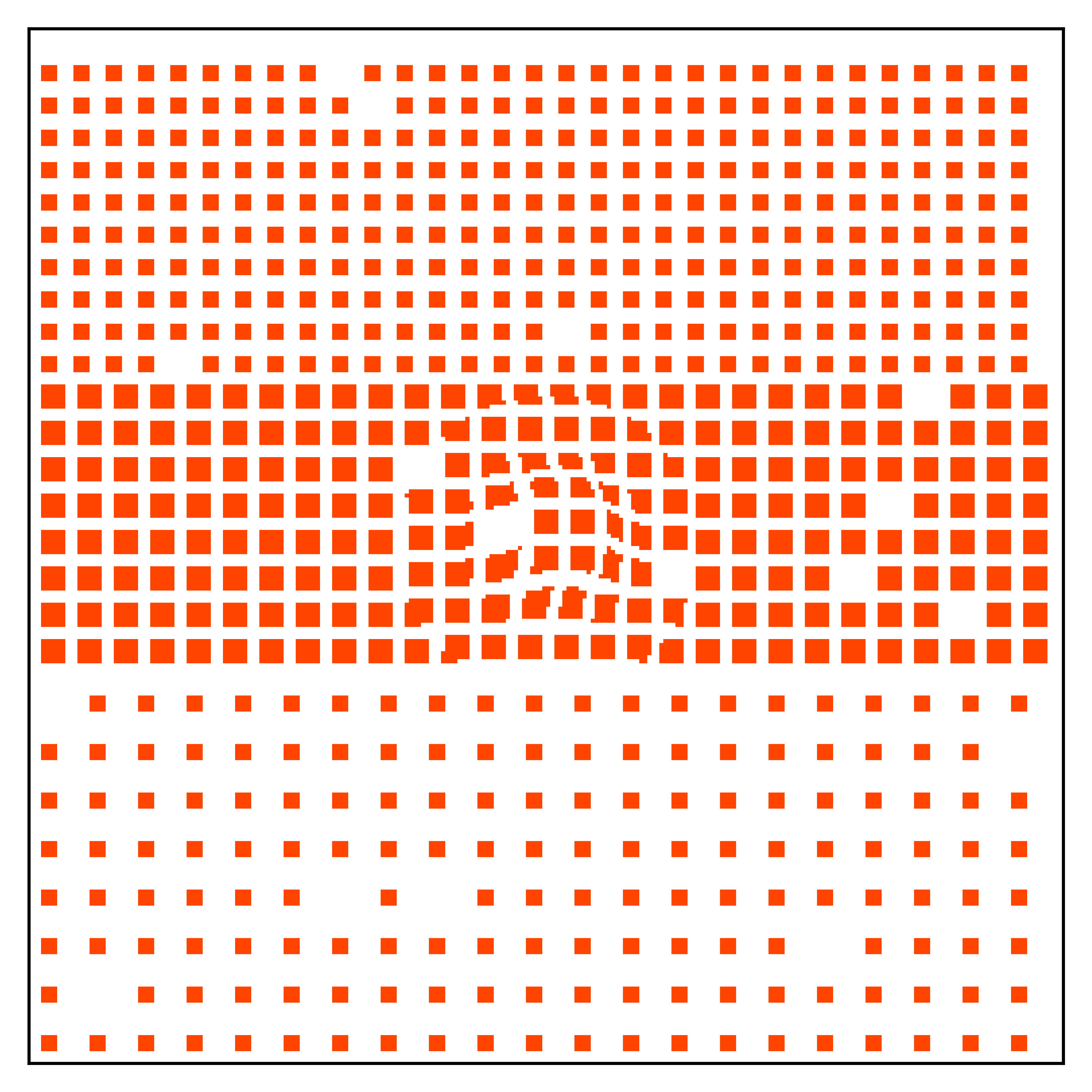}
	\end{subfigure}
	\begin{subfigure}[b]{0.375\textwidth}
				\hbox{\vspace{-0.365em}
		\includegraphics[width=\textwidth]{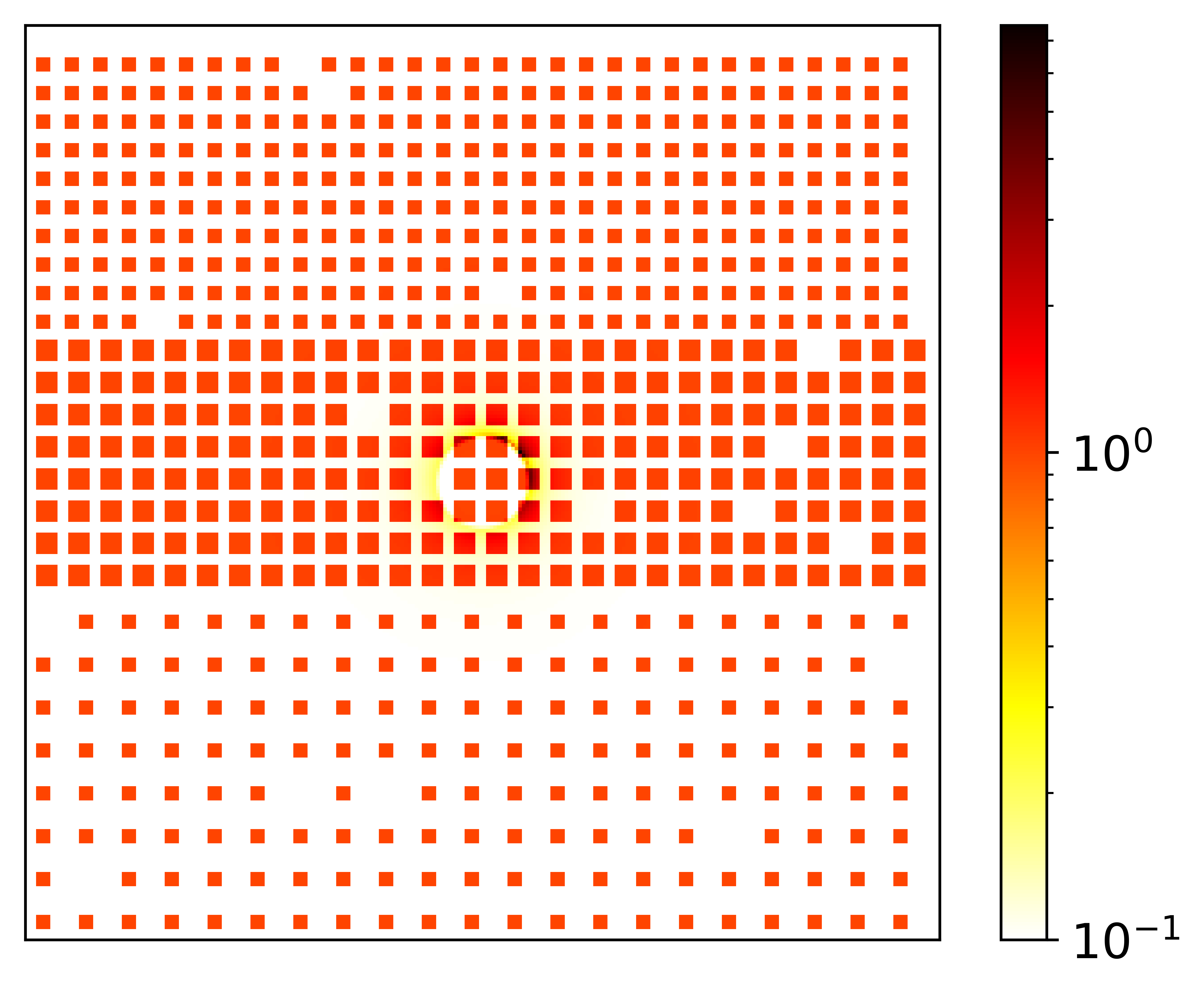}}
	\end{subfigure}
	\centering
	
	\caption{Reference coefficient $A_{\text{ref}}$ (left), perturbation in the physical domain $A_y$ (center) and corresponding change in value perturbation $A$ (right).}
	\label{ex2:coefficients}
\end{figure}

\begin{figure}[h]
	\begin{subfigure}[b]{0.32\textwidth}
		\includegraphics[width=\textwidth]{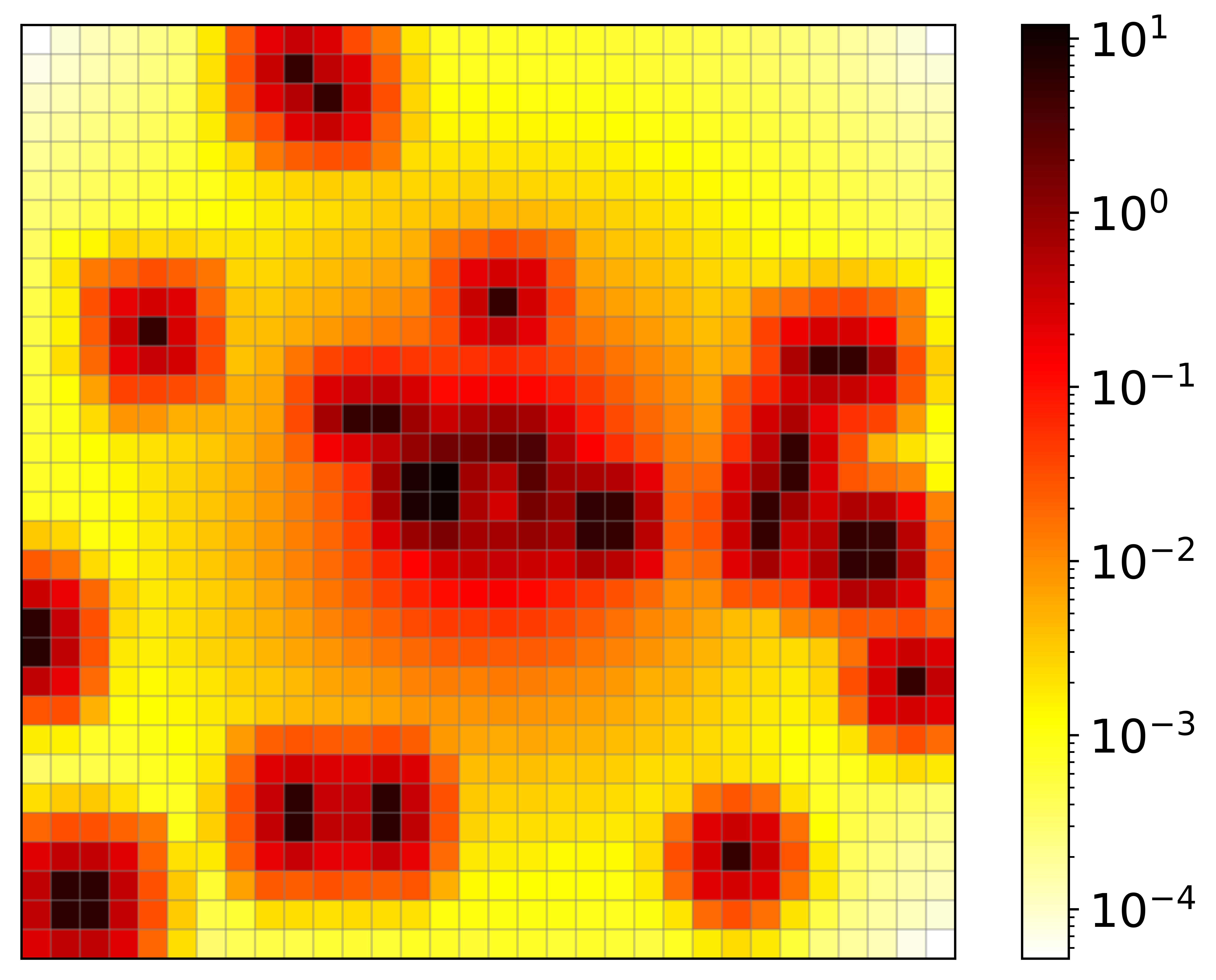}
	\end{subfigure}
	\begin{subfigure}[b]{0.32\textwidth}
		\includegraphics[width=\textwidth]{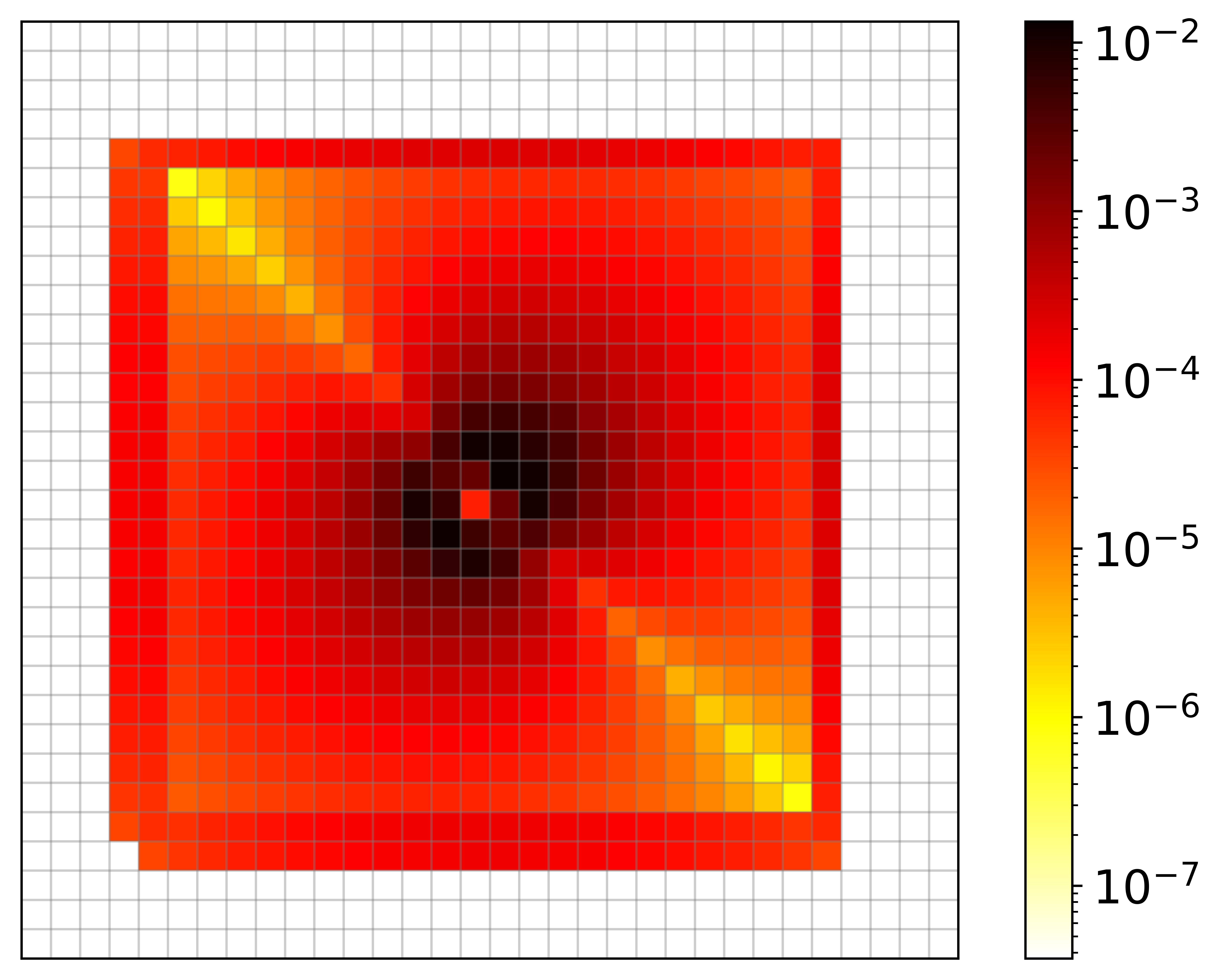}
	\end{subfigure}
	\begin{subfigure}[b]{0.32\textwidth}
		\includegraphics[width=\textwidth]{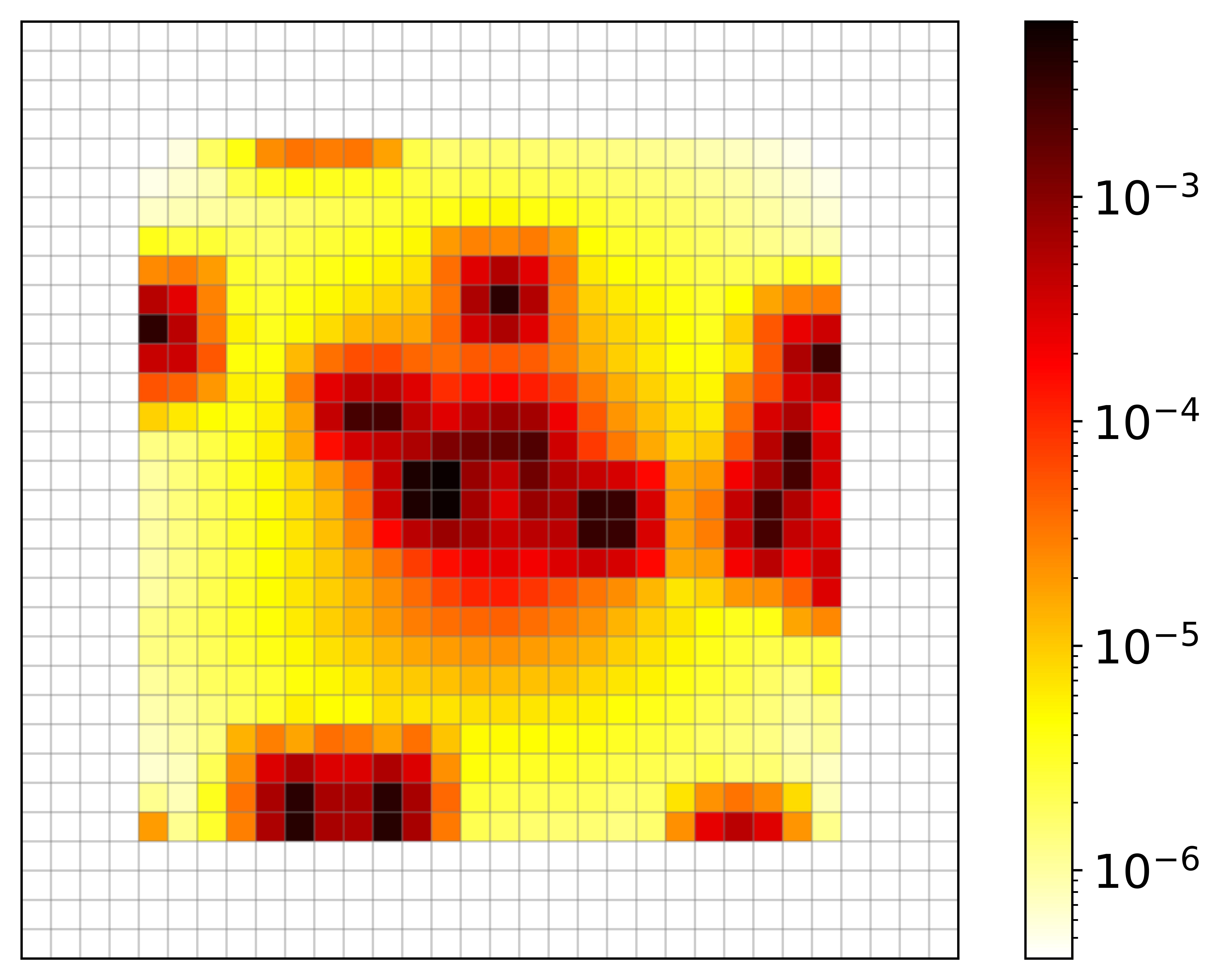}
	\end{subfigure}
	\centering
	\caption{Error indicators $E_{\QQ V_H,T} \cdot \norm{f}_{L^2(\Omega)}$ (left), ${E_{f,T}}$  (center) and ${E_{\RR f,T}}$ (right).}
	\label{ex2:indicators}
\end{figure}

\begin{figure}[h]
	\centering
	\includegraphics[width=0.55\textwidth]{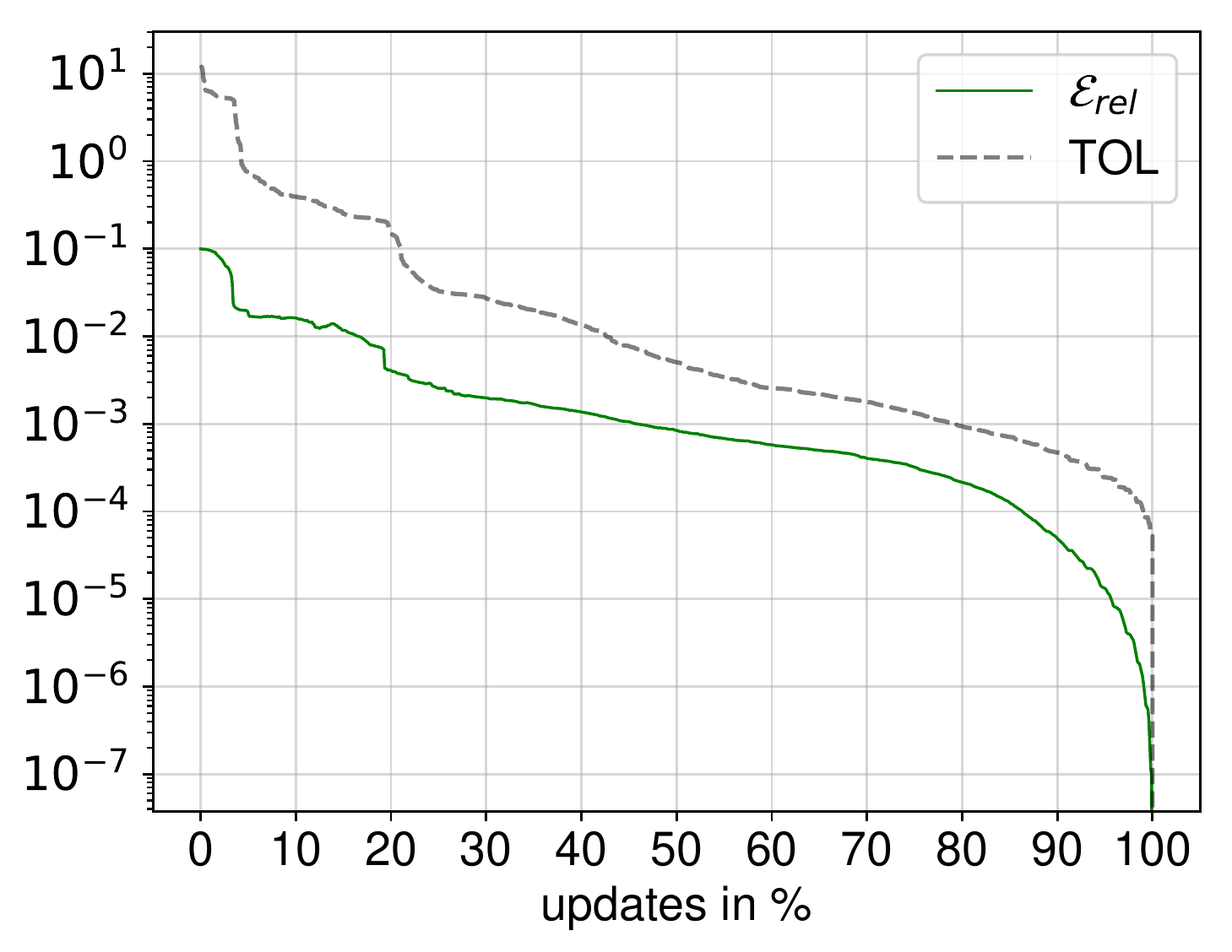}
	\caption{Relative error improvement for local domain mapping.}
	\label{ex2:errors}
\end{figure}

In the second experiment we choose $\psi$ to be a local distortion in the middle of the domain which can be seen in Figure \ref{ex2:coefficients}. 
With the help of domain mappings the reference coefficient $A_{\text{ref}}$ is subjected to a simple change in value which means that $\Omega_1$ does not change its position. 
This is visualized in the right picture of Figure \ref{ex2:coefficients}. 
The domain mapping as well as the defects can be clearly seen in $E_{\QQ V_H,T}$, whereas the defects stick out compared to the domain mapping. 
In contrast to the diffusion coefficient, $f$ is a distortion of $f_{\text{ref}} = f_y$. 
This means that the support of the two is not the same which can be seen in ${E_{f,T}}$ in Figure \ref{ex2:indicators} . Furthermore $f$ is not affected by the defects in $A$ which explains that ${E_{f,T}}$ only detects the domain mapping and not the defects.
In a view of Figure \ref{ex2:errors}, we observe a similar effect as in the first experiment. 
However, due to the domain mapping, it takes comparatively longer to converge to the optimal PGLOD solution. 

\clearpage
\subsection{Global domain mapping}
\begin{figure}[h]
	\begin{subfigure}[b]{0.3\textwidth}
		\includegraphics[width=\textwidth]{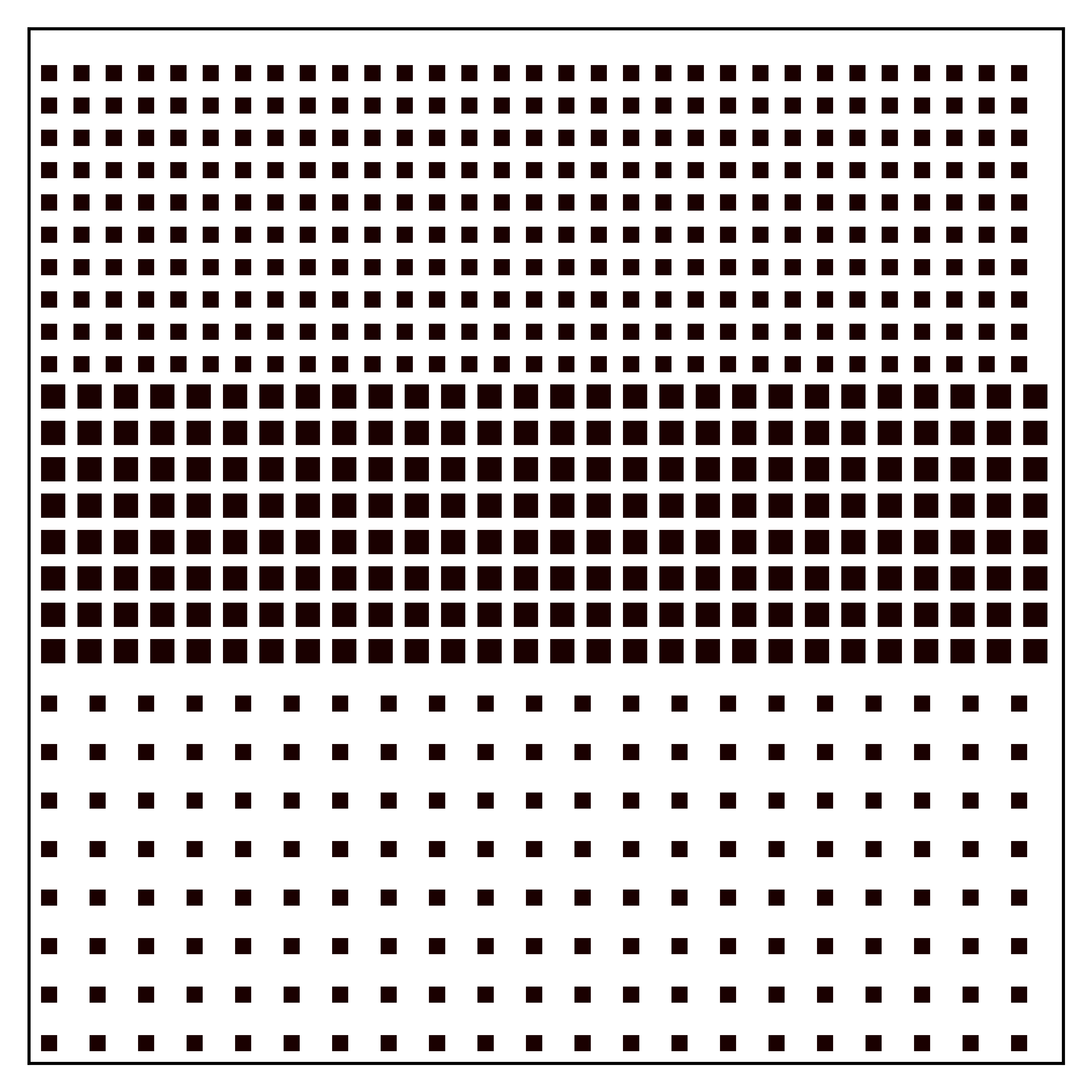}
	\end{subfigure}
	\begin{subfigure}[b]{0.3\textwidth}
		\includegraphics[width=\textwidth]{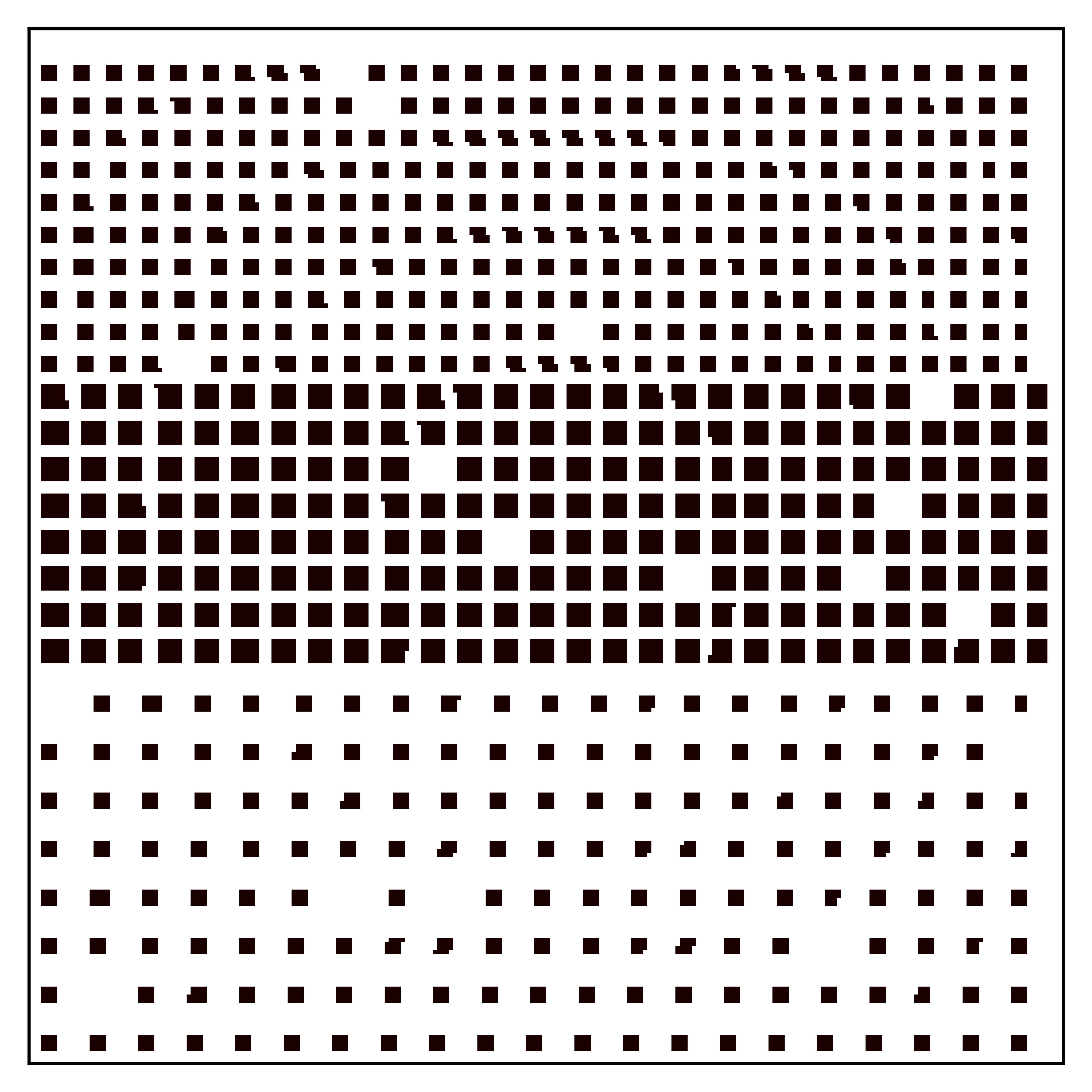}
	\end{subfigure}
	\begin{subfigure}[b]{0.375\textwidth}
			\hbox{\vspace{-0.365em}
	\includegraphics[width=\textwidth]{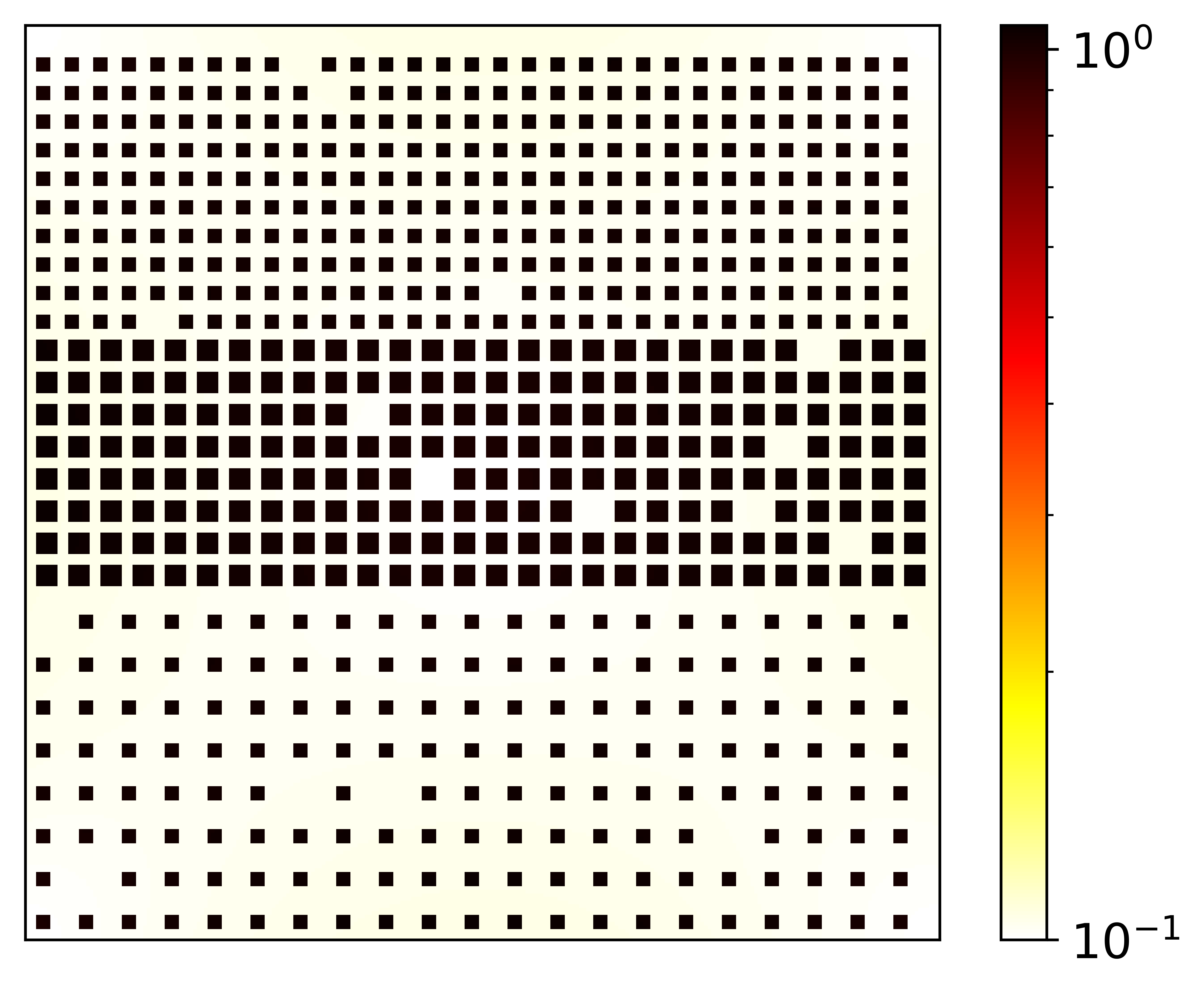}}
	\end{subfigure}
	\centering
	\caption{Reference coefficient $A_{\text{ref}}$ (left), perturbation in the physical domain $A_y$ (center) and corresponding change in value perturbation $A$ (right).}
	\label{ex3:coeficients}	
\end{figure}

\begin{figure}[h]
	\begin{subfigure}[b]{0.32\textwidth}
		\includegraphics[width=\textwidth]{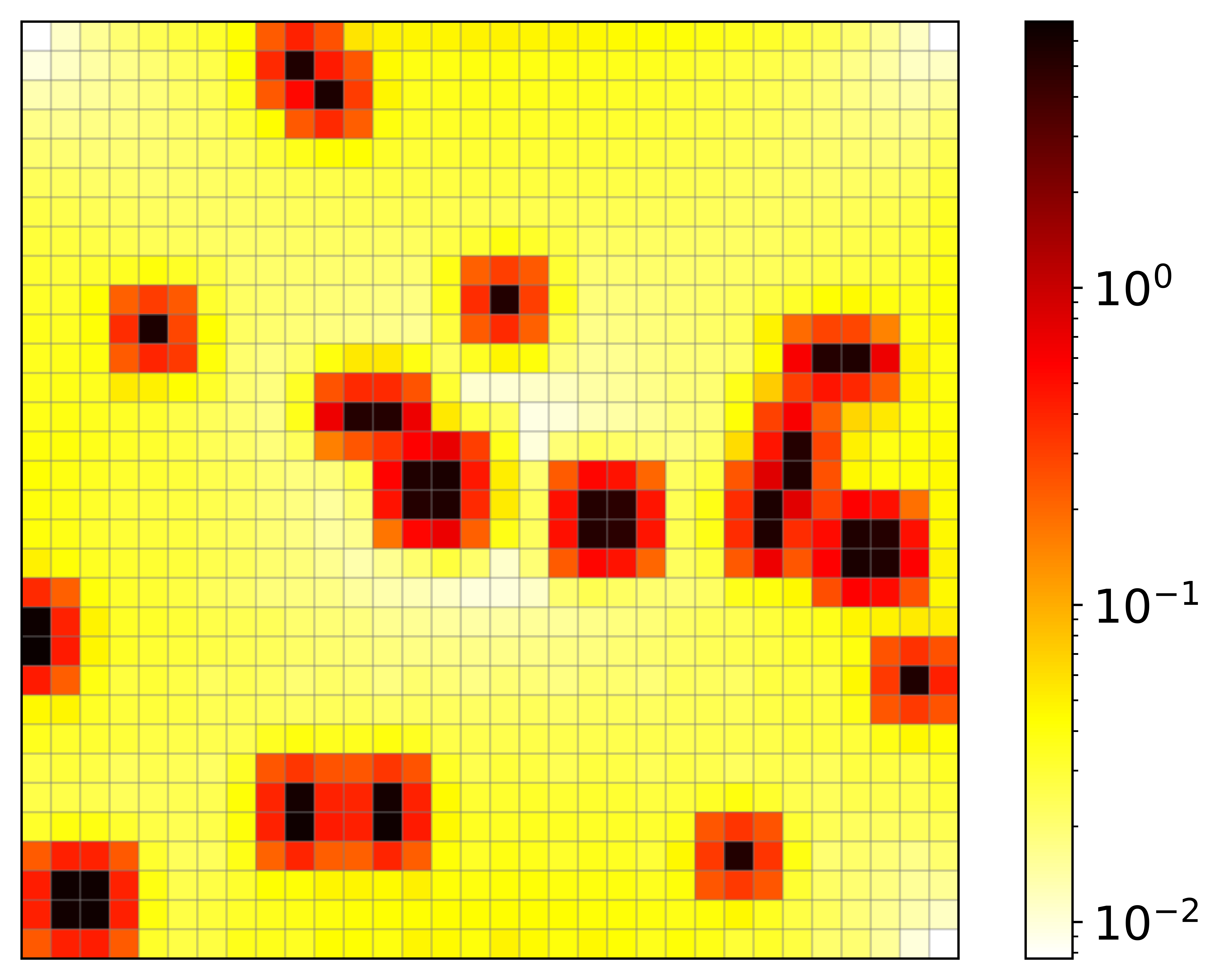}
	\end{subfigure}
	\begin{subfigure}[b]{0.32\textwidth}
		\includegraphics[width=\textwidth]{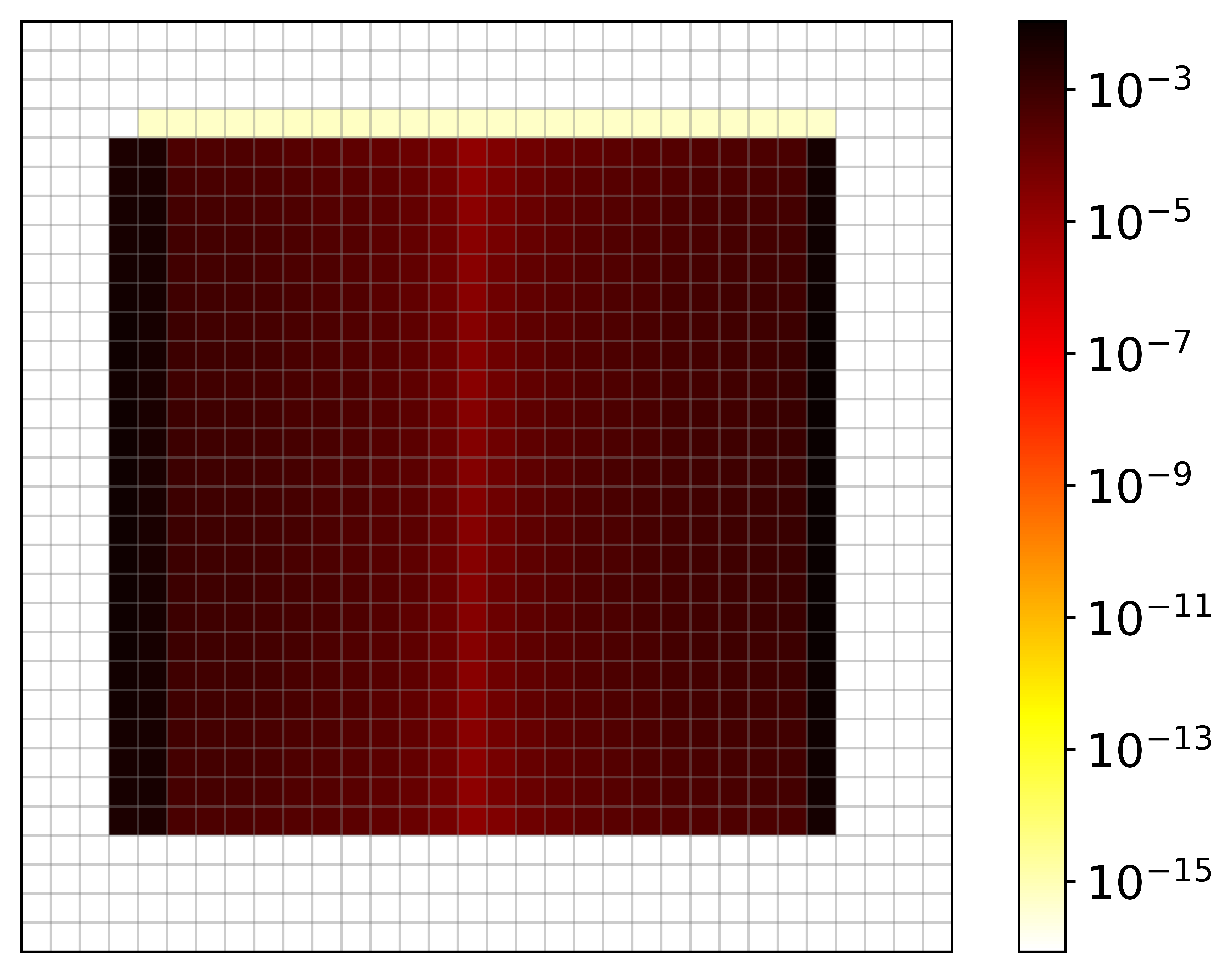}
	\end{subfigure}
	\begin{subfigure}[b]{0.32\textwidth}
		\includegraphics[width=\textwidth]{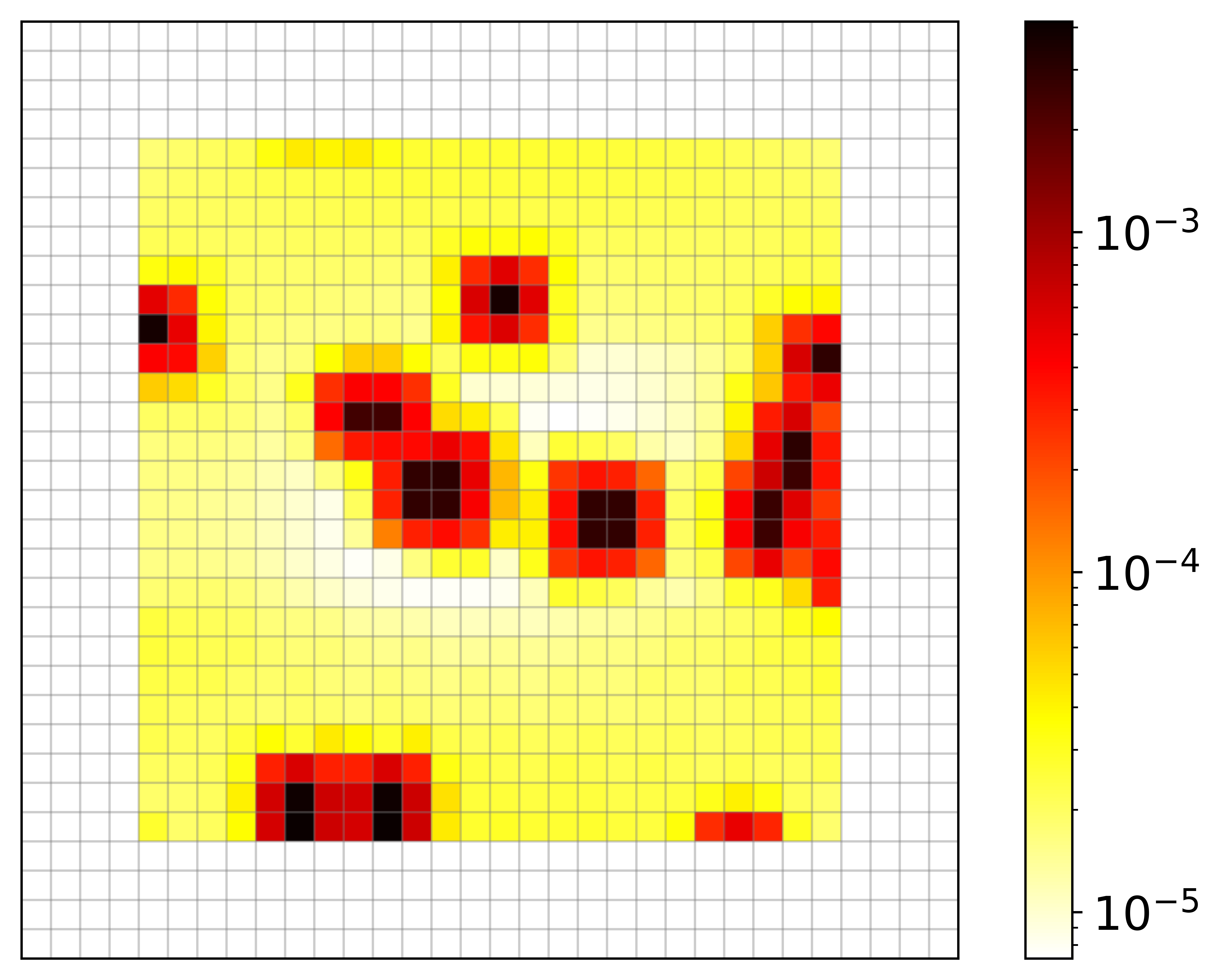}
	\end{subfigure}
	\centering
	\caption{Error indicators $E_{\QQ V_H,T} \cdot \norm{f}_{L^2(\Omega)}$ (left), ${E_{f,T}}$  (center) and ${E_{\RR f,T}}$ (right).}
		\label{ex3:indicators}
\end{figure}

In our third experiment we address the situation that $\psi$ has support on the whole domain $\Omega$. 
The perturbation in the physical domain as well as the change in value of the reference coefficient can be seen in Figure \ref{ex3:coeficients}. 
The yellow background in $E_{\QQ V_H,T}$ in Figure \ref{ex3:indicators} visualizes the effect of the domain mapping and the defects are clearly notable. 
In ${E_{f,T}}$ we once again see the change of support in $f$ (left and right black channel) and the defects are visible in ${E_{\RR f,T}}$. 

\begin{figure}[h]
	\centering
	\includegraphics[width=0.55\textwidth]{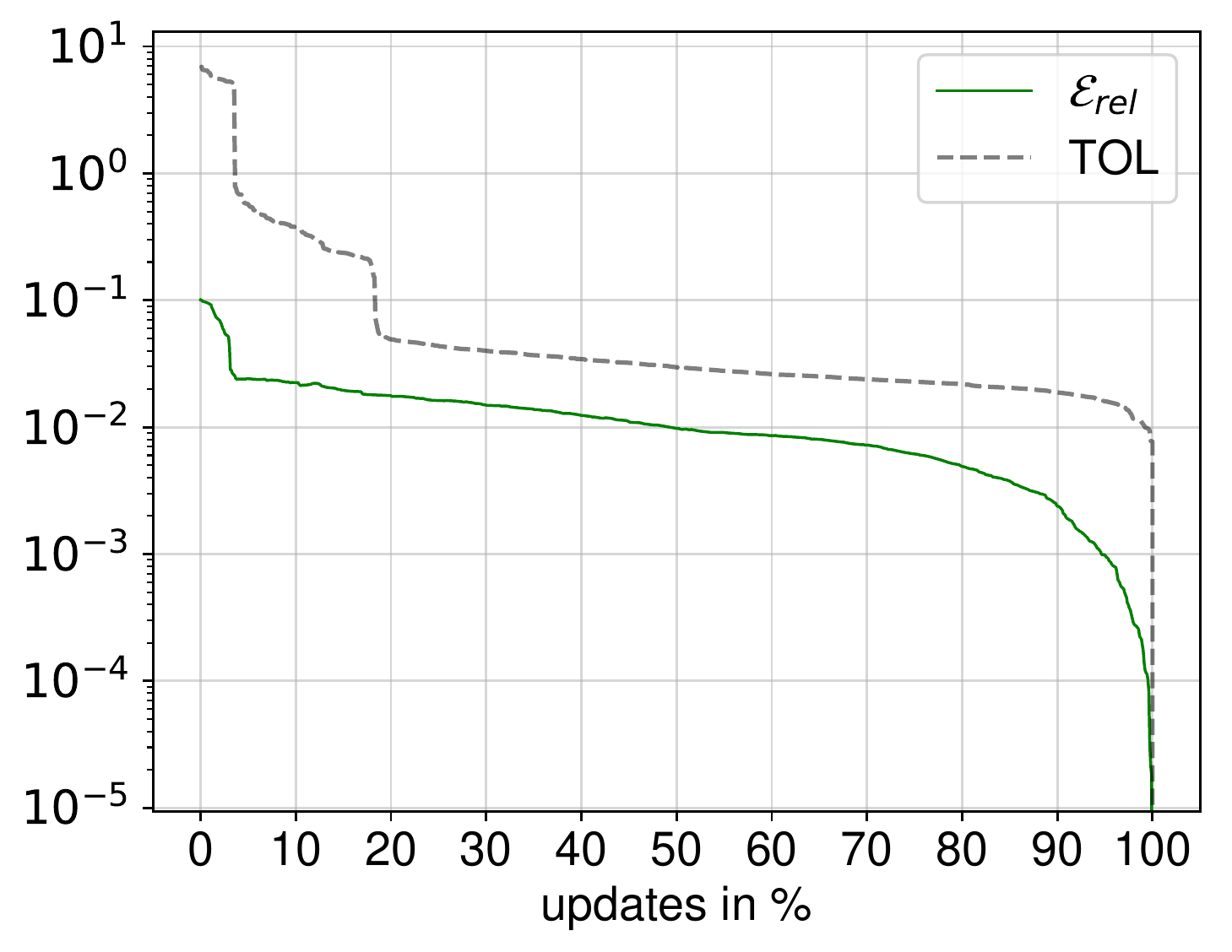}
	\caption{Relative error improvement for global domain mapping.}
	\label{ex3:errors}
\end{figure}

In Figure \ref{ex3:errors}, we observe that local defects are still resolved efficiently whereas the global map causes a relatively low convergence to the optimal PGLOD solution. 
Thus, this example shows that global domain mappings are difficult to handle. 
However, for instance for the case that an accuracy of $10^{-2}$ is accurate enough, we still get a reasonable result. 
Note that in this particular example no use of domain mappings would surely result in $100\%$ recomputation as the complete coefficient changes.   

% %\axel{I think we need something like this but not with 1. and 2.  perhaps. Its the same problem for standard FEM when we have the bound $e\leq CH$ and we pick an H but we kow the constant is too high so its hard to say what error a certain H corresponds to. Some kind of comparison  between different approximation levels is needed. Perhaps we should express it a bit more vaguely without TOL/2.}
%\tim{I changed this now}
\newpage
\begin{bemerkung}[Choosing TOL]
	In all our experiments, $\mathcal{E}_{\text{rel}}(u_k,\tilde{u}_k)$ shows a promising behavior. 
	However, we point out that in practice, we clearly do not know $\mathcal{E}_{\text{rel}}(u_k,\tilde{u}_k)$ a priori. 
	Thus, we start with a rather high tolerance TOL and consider an update error $\mathcal{E}_{\text{rel}}(\tilde{u}_{k,\text{old}},\tilde{u}_k)$, 
	where $\tilde{u}_{k,\text{old}}$ denotes on old approximation with respect to the former tolerance. 
	Whenever the number of updates from one tolerance to another is strictly positive and $\mathcal{E}_{\text{rel}}(\tilde{u}_{k,\text{old}},\tilde{u}_k)$ does still exhibit a significant gain, 
	we should proceed with a smaller TOL and continue this algorithm until $\mathcal{E}_{\text{rel}}(\tilde{u}_{k,\text{old}},\tilde{u}_k)$ is small enough. 
\end{bemerkung}

\bibliographystyle{plain}

\end{document}